\documentclass[11pt]{amsart}

\usepackage[margin=1.1in]{geometry}
\usepackage{amsmath,amssymb,amsthm,mathtools}
\usepackage{enumitem}
\usepackage{booktabs}
\usepackage{xcolor}
\usepackage{microtype}
\usepackage[most]{tcolorbox}
\definecolor{linkblue}{RGB}{20,60,130}
\definecolor{boxgray}{RGB}{246,247,250}
\definecolor{rulegray}{RGB}{120,135,165}
\usepackage[colorlinks=true,linkcolor=linkblue,citecolor=linkblue,urlcolor=linkblue]{hyperref}

\newtheorem{theorem}{Theorem}[section]
\newtheorem{lemma}[theorem]{Lemma}
\newtheorem{proposition}[theorem]{Proposition}
\newtheorem{corollary}[theorem]{Corollary}
\newtheorem{mainthm}{Theorem}

\tcolorboxenvironment{mainthm}{enhanced,breakable,colback=boxgray,colframe=boxgray,
  borderline west={1.2pt}{0pt}{rulegray},arc=0pt,outer arc=0pt,
  left=8pt,right=6pt,top=2pt,bottom=4pt,before skip=10pt,after skip=10pt}
\theoremstyle{definition}
\newtheorem{definition}[theorem]{Definition}
\newtheorem{example}[theorem]{Example}
\theoremstyle{remark}
\newtheorem{remark}[theorem]{Remark}

\numberwithin{equation}{section}

\setlist{itemsep=2pt,topsep=4pt}
\allowdisplaybreaks

\newcommand{\N}{\mathbb{N}}
\newcommand{\R}{\mathbb{R}}
\newcommand{\C}{\mathbb{C}}
\newcommand{\D}{\mathbb{D}}
\newcommand{\Fix}{\operatorname{Fix}}
\newcommand{\sgn}{\operatorname{sgn}}
\newcommand{\BMOA}{\mathrm{BMOA}}
\newcommand{\norm}[1]{\left\lVert #1\right\rVert}
\newcommand{\abs}[1]{\left\lvert #1\right\rvert}

\title[Rational contractions in quasi-Banach spaces]{Rational contractions in quasi-Banach spaces: convergence, multiplicity of fixed points and applications}

\author[A. H. Ansari Komachali, Y. U. Gaba, M. Aphane, I. Yildirim]{Arslan Hojat Ansari Komachali$^{\mathrm{a}}$, Ya\'e Ulrich Gaba$^{\mathrm{a},\mathrm{b},\mathrm{c},*}$, Maggie Aphane$^{\mathrm{a}}$, Isa Yildirim$^{\mathrm{d}}$}

\dedicatory{\upshape\small
$^{\mathrm{a}}$Department of Mathematics and Applied Mathematics, Sefako Makgatho Health Sciences University, Ga-Rankuwa, Pretoria, Medunsa-0204, South Africa\\[2pt]
$^{\mathrm{b}}$AI Research and Innovation Nexus for Africa (AIRINA Labs), B\'enin\\[2pt]
$^{\mathrm{c}}$African Center for Advanced Studies (ACAS), P.O. Box 4477, Yaound\'e, Cameroon\\[2pt]
$^{\mathrm{d}}$Department of Mathematics, Faculty of Science, Atat\"urk University, Erzurum 25240, Turkey}

\thanks{$^{*}$Corresponding author: Ya\'e Ulrich Gaba (\href{mailto:yaeulrich.gaba@gmail.com}{yaeulrich.gaba@gmail.com}).}

\subjclass[2020]{Primary 47H10, 54H25; Secondary 46A16, 47G10, 30H10, 39B82}
\keywords{Fixed point; $b$-metric space; quasi-Banach space; rational contraction; convergence rate; multiplicity of fixed points; Ulam--Hyers stability; $\ell^p$ with $0<p<1$; Hardy space $H^p$; Volterra-type integration operator}

\date{}

\makeatletter
\def\@setauthors{%
  \begingroup
  \def\thanks{\protect\thanks@warning}%
  \trivlist
  \centering\footnotesize \@topsep30\p@\relax
  \advance\@topsep by -\baselineskip
  \item\relax
  \normalsize\authors
  \endtrivlist
  \endgroup}
\makeatother

\begin{document}

\begin{abstract}
We study self-maps $f$ of a complete $b$-metric space $(U,d)$ (constant $b\ge1$) satisfying the rational contractive condition
\[
d(fu,fv)\le\frac{d(u,fv)+d(v,fu)+\kappa}{b\,[d(u,fu)+d(v,fv)]+1}\,d(u,v)\qquad(u,v\in U),
\]
with $\kappa\in[0,1)$. We first observe that the two-parameter family of conditions studied in our earlier work reduces, after a rescaling of the $b$-metric, to this one-parameter form. We then prove: (A) every Picard sequence converges to a fixed point, with an explicit a~priori estimate whose geometric ratio depends only on $b$, $\kappa$ and $d(u_0,fu_0)$, an a~posteriori estimate, and asymptotic rate at most $\kappa$; (B) a structure theorem: such maps are automatically continuous, every fixed point attracts the whole ball of radius $(1-\kappa)/2$, distinct fixed points are at distance at least $(1-\kappa)/2$ (sharp), basins of attraction are clopen, and the number of fixed points does not exceed the number of connected components of $U$; in particular the fixed point is unique on connected spaces; (C) Ulam--Hyers stability of the fixed-point problem. The class of such maps is not invariant under passing from a quasi-norm to the equivalent $p$-norm metric, so the quasi-Banach results are not formal consequences of the metric ones. A gluing principle shows that the multiplicity permitted by (B) occurs. As applications we obtain existence, exact multiplicity and error estimates for infinite nonlinear systems with relay-type forcing in $\ell^p$, $0<p<1$, and for nonlinear perturbations of Volterra-type equations $f=\pm a+T_gf+\psi(f(z_0))k$ in the Hardy space $H^p$, $0<p<1$, where $T_gf(z)=\int_0^z f(\zeta)g'(\zeta)\,d\zeta$ and $g\in\BMOA$.
\end{abstract}

\maketitle

\section{Introduction}

Since Bakhtin~\cite{Bakhtin} and Czerwik~\cite{Czerwik93,Czerwik98}, $b$-metric spaces have become a standard setting for metric fixed point theory. Two facts limit what can be new in this setting. First, the inequality of Lemma~\ref{lem:geom} below (in the form proved by Miculescu--Mihail~\cite{MM}, Suzuki~\cite{Suzuki} and Huang--Deng--Radenovi\'c~\cite{HDR}) shows that a Picard sequence with $d(u_n,u_{n+1})\le\delta\,d(u_{n-1},u_n)$ is Cauchy for every $\delta<1$, independently of $b$. Second, by the metrization theorem of Paluszy\'nski--Stempak~\cite{PS} (going back to Aoki~\cite{Aoki} and Rolewicz~\cite{Rolewicz} for quasi-norms), every $b$-metric is equivalent to a power of a metric. Consequently many fixed point statements in $b$-metric spaces are formal consequences of their metric counterparts; see the discussion in Kirk--Shahzad~\cite{KS}.

This paper continues the study of \emph{rational-type contractions}, in which the contraction factor is a ratio of displacements. The idea goes back to Khojasteh, Abbas and Costache~\cite{KAC}, who worked in metric spaces (Theorem~\ref{thm:KAC}); Demma and Vetro~\cite{DV} carried it over to $b$-metric spaces (Theorem~\ref{thm:DV}), and in~\cite{AGAY} we relaxed it further by allowing maps that satisfy
\begin{equation}\label{eq:old}
d(fu,fv)\le\frac{d(u,fv)+d(v,fu)+p}{b\,[d(u,fu)+d(v,fv)]+m}\,d(u,v),\qquad 0<p<m
\end{equation}
(Theorems~\ref{thm:old22} and~\ref{thm:old23}). These four results, recalled in Section~\ref{subsec:known}, have the same conclusion: a fixed point exists, Picard iterates converge to a fixed point, and distinct fixed points are separated by a fixed distance. They do not give the rate of convergence, nor do they say when several fixed points can occur. We address both questions. Our main examples are the quasi-Banach spaces $\ell^p$ and $H^p$ with $0<p<1$, where the quasi-metric comes from the space itself. We do not use $L^p(0,1)$, $0<p<1$: its elements need not be locally integrable and its dual is trivial~\cite{Day}, so integral operators are of little use there.

We first observe that the two constants in~\eqref{eq:old} are not independent: dividing the $b$-metric by $m$ turns~\eqref{eq:old} into a condition with the single parameter $\kappa=p/m$ (Proposition~\ref{prop:reduction}). This one-parameter condition $(\mathrm R_\kappa)$ of Definition~\ref{def:rational} is studied throughout. Section~\ref{sec:A} deals with convergence: every Picard sequence converges with a geometric ratio that depends only on $b$, $\kappa$ and the initial displacement $d(u_0,fu_0)$, and we obtain an a priori bound, an a posteriori bound, and an asymptotic rate at most $\kappa$ (Theorem~\ref{thm:A}), none of which follows from Theorems~\ref{thm:KAC}--\ref{thm:old23}. Example~\ref{ex:nonBanach} gives a map on $[0,\infty)$ to which Theorem~\ref{thm:A} applies but the theorems of Banach, Boyd--Wong and Meir--Keeler do not. Section~\ref{sec:structure} turns to part~(iii) of those theorems, which allows several fixed points without saying where they lie. We show that $(\mathrm R_\kappa)$ forces the map to be continuous, every fixed point attracts the whole ball of radius $(1-\kappa)/2$, and the basins of attraction are open and closed, so that the number of fixed points cannot exceed the number of connected components of the space (Theorem~\ref{thm:B}). On a connected space, for instance a closed convex subset of a quasi-Banach space, the fixed point is unique, and the same estimates make the fixed-point problem Ulam--Hyers stable (Theorem~\ref{thm:C}, Section~\ref{sec:C}). Section~\ref{sec:consequences} draws the consequences for Theorems~\ref{thm:KAC}--\ref{thm:old23}, and shows by an example that no such structure result can hold for the variants with an additional term $L\,d(v,fu)$, $L\ge1$. In Section~\ref{sec:choice} we check that these results in quasi-Banach spaces do not follow from their metric versions through the Aoki--Rolewicz theorem: since $(\mathrm R_\kappa)$ is not homogeneous, the class of maps it defines changes when a quasi-norm is replaced by the equivalent $p$-norm metric. The multiplicity allowed by Theorem~\ref{thm:B} does occur: in Section~\ref{sec:gluing}, a gluing principle (Theorem~\ref{thm:D}) assembles rational contractions with a prescribed number of fixed points out of uniformly contractive pieces, and we use it to obtain exactly two solutions, with explicit error estimates, for an infinite nonlinear system in $\ell^p$ (Theorem~\ref{thm:E}, Section~\ref{sec:lp}) and for a nonlinear perturbation of a Volterra-type equation in the Hardy space $H^p$ (Theorem~\ref{thm:F}, Section~\ref{sec:Hp}), both with $0<p<1$. Section~\ref{sec:open} contains concluding remarks.

\section{Preliminaries}

We recall $b$-metric spaces and the fact that their topology is metrizable, which is needed to speak of continuity, connectedness and compactness. We then prove a quantitative version of the Cauchy criterion for Picard sequences, introduce quasi-Banach and $p$-Banach spaces, state precisely the results we extend, and finally define the class of maps studied in the paper.

\subsection{\texorpdfstring{$b$}{b}-metric spaces and their topology}

\begin{definition}[\cite{Czerwik98}]
Let $U\neq\emptyset$ and $b\ge1$. A map $d:U\times U\to[0,\infty)$ is a \emph{$b$-metric} if, for all $u,v,y\in U$: $d(u,v)=0\iff u=v$; $d(u,v)=d(v,u)$; and $d(u,y)\le b\,[d(u,v)+d(v,y)]$.
\end{definition}

A $b$-metric with constant $b$ is also a $b$-metric with every constant $b'\ge b$; in what follows $b$ is part of the data. For $b=1$ we recover metric spaces.

\begin{definition}[\cite{Czerwik98}]\label{def:conv}
Let $(U,d)$ be a $b$-metric space and $(u_n)\subseteq U$.
\begin{enumerate}[label=(\alph*)]
\item $(u_n)$ \emph{converges} to $u\in U$, written $u_n\to u$, if $d(u_n,u)\to0$;
\item $(u_n)$ is \emph{Cauchy} if $d(u_n,u_m)\to0$ as $n,m\to\infty$;
\item $(U,d)$ is \emph{complete} if every Cauchy sequence converges.
\end{enumerate}
\end{definition}

A $b$-metric need not be continuous, and open balls need not be open sets. We give $U$ the topology $\tau_d$ in which $G$ is open iff for every $x\in G$ there is $\varepsilon>0$ with $\left\lbrace y:d(x,y)<\varepsilon\right\rbrace \subseteq G$.

\begin{lemma}[Paluszy\'nski--Stempak~\cite{PS}]\label{lem:PS}
Let $(U,d)$ be a $b$-metric space and let $r\in(0,1]$ be defined by $(2b)^r=2$. Then
\[
\rho(x,y):=\inf\left\lbrace \sum_{i=0}^{n}d(x_i,x_{i+1})^r:\ n\in\N,\ x_0=x,\ x_{n+1}=y\right\rbrace 
\]
is a metric on $U$ satisfying $\rho(x,y)\le d(x,y)^r\le 4\rho(x,y)$ for all $x,y\in U$.
\end{lemma}

\begin{lemma}\label{lem:topology}
With $\rho$ as in Lemma~\ref{lem:PS}:
\begin{enumerate}[label=(\alph*)]
\item $\tau_d$ coincides with the topology of $\rho$; in particular it is metrizable, and continuity with respect to $\tau_d$ is equivalent to sequential continuity;
\item $d(u_n,u)\to0$ iff $\rho(u_n,u)\to0$; limits are unique; $(u_n)$ is $d$-Cauchy iff it is $\rho$-Cauchy; and $(U,d)$ is complete iff $(U,\rho)$ is complete.
\end{enumerate}
\end{lemma}

\begin{proof}
By Lemma~\ref{lem:PS}, $\left\lbrace y:\rho(x,y)<\varepsilon^r/4\right\rbrace \subseteq\left\lbrace y:d(x,y)<\varepsilon\right\rbrace \subseteq\left\lbrace y:\rho(x,y)<\varepsilon^r\right\rbrace $. This gives (a), and (b) follows from the same two-sided inequality.
\end{proof}

\begin{remark}\label{rem:pnorm}
If $d^{\,r}$ itself satisfies the triangle inequality, one may take $\rho=d^{\,r}$ in Lemma~\ref{lem:PS}, and then the constant $4$ may be replaced by $1$ in everything that follows. This happens for $p$-normed spaces (Section~\ref{subsec:quasi}), where $b=2^{1/p-1}$ and $r=p$.
\end{remark}

\subsection{A geometric Cauchy estimate}

For a long time the standard tool for Picard sequences in $b$-metric spaces was the criterion of Jovanovi\'c, Kadelburg and Radenovi\'c~\cite{JKR}, which requires $\delta<1/b$. It was later shown that the restriction is unnecessary.

\begin{lemma}[Miculescu--Mihail~\cite{MM}, Suzuki~\cite{Suzuki}, Huang--Deng--Radenovi\'c~\cite{HDR}]\label{lem:known}
Let $(U,d)$ be a $b$-metric space and $(u_n)\subseteq U$. If there is $\delta\in[0,1)$ with $d(u_n,u_{n+1})\le\delta\,d(u_{n-1},u_n)$ for all $n\ge1$, then $(u_n)$ is Cauchy.
\end{lemma}

The next lemma refines Lemma~\ref{lem:known} by giving an explicit bound; we include its short proof.

\begin{lemma}\label{lem:geom}
Let $(U,d)$ be a $b$-metric space, let $r$ be as in Lemma~\ref{lem:PS}, and set
\[
C_\theta:=\left(\frac{4}{1-\theta^r}\right)^{1/r},\qquad\theta\in[0,1).
\]
If $(u_n)_{n\ge0}\subseteq U$ satisfies $d(u_n,u_{n+1})\le c\,\theta^n$ for all $n\ge0$, where $c\ge0$, then $d(u_n,u_m)\le C_\theta\,c\,\theta^n$ for all $m>n\ge0$. In particular $(u_n)$ is Cauchy, and if $u_n\to u$ then $d(u_n,u)\le C_\theta\,c\,\theta^n$ for all $n$.
\end{lemma}

\begin{proof}
By the triangle inequality for $\rho$ and Lemma~\ref{lem:PS},
\[
\rho(u_n,u_m)\le\sum_{k=n}^{m-1}\rho(u_k,u_{k+1})\le\sum_{k=n}^{m-1}d(u_k,u_{k+1})^r\le c^r\sum_{k\ge n}\theta^{kr}=\frac{c^r\theta^{nr}}{1-\theta^r},
\]
and $d(u_n,u_m)^r\le4\rho(u_n,u_m)$. If $u_m\to u$, then $\rho(u_m,u)\le d(u_m,u)^r\to0$, so letting $m\to\infty$ in $\rho(u_n,u)\le\rho(u_n,u_m)+\rho(u_m,u)$ gives $\rho(u_n,u)\le c^r\theta^{nr}/(1-\theta^r)$, and again $d(u_n,u)^r\le4\rho(u_n,u)$.
\end{proof}

\subsection{Quasi-Banach and \texorpdfstring{$p$}{p}-Banach spaces}\label{subsec:quasi}

A quasi-norm on a vector space $X$ is a map $\norm{\cdot}$ satisfying the norm axioms with the triangle inequality replaced by $\norm{x+y}\le K(\norm x+\norm y)$ for some $K\ge1$ (the \emph{modulus}). Then $d(x,y)=\norm{x-y}$ is a $b$-metric with $b=K$, and $X$ is \emph{quasi-Banach} if it is complete; see~\cite{KPR}. For $0<p\le1$, $X$ is \emph{$p$-normed} if $\norm{x+y}^p\le\norm x^p+\norm y^p$. In that case $\norm{\cdot}^p$ is a translation-invariant metric, and since $(a^p+b^p)^{1/p}\le2^{1/p-1}(a+b)$ for $a,b\ge0$, the modulus is $K=2^{1/p-1}$. With $b=K$ we have $(2b)^p=2$, so $r=p$ and Remark~\ref{rem:pnorm} applies. The spaces $\ell^p$ and $H^p$ with $0<p<1$ are $p$-Banach spaces. In a $p$-normed space we use the two elementary consequences
\begin{equation}\label{eq:pineq}
\norm{x+y}^p\le\norm x^p+\norm y^p,\qquad\norm{x+y}^p\ge\norm x^p-\norm y^p.
\end{equation}

\subsection{Known results for rational-type contractions}\label{subsec:known}

We now state the results that the present paper extends. In all four, the conclusion has the same form:
\begin{enumerate}[label=(\roman*)]
\item $f$ has at least one fixed point;
\item every Picard sequence $u_n=f^nu_0$, $u_0\in U$, converges to a fixed point of $f$;
\item distinct fixed points $u^*\ne v^*$ satisfy $d(u^*,v^*)\ge\Delta$,
\end{enumerate}
where the separation constant $\Delta$ depends on the theorem.

\begin{theorem}[Khojasteh--Abbas--Costache~\cite{KAC}]\label{thm:KAC}
Let $(U,d)$ be a complete metric space and let $f:U\to U$ satisfy
\[
d(fu,fv)\le\frac{d(u,fv)+d(v,fu)}{d(u,fu)+d(v,fv)+1}\,d(u,v)\qquad(u,v\in U).
\]
Then \textup{(i)--(iii)} hold with $\Delta=1/2$.
\end{theorem}

\begin{theorem}[Demma--Vetro~\cite{DV}]\label{thm:DV}
Let $(U,d)$ be a complete $b$-metric space, $L\ge0$, and let $f:U\to U$ satisfy
\[
d(fu,fv)\le\frac{d(u,fv)+d(v,fu)}{b[d(u,fu)+d(v,fv)]+1}\,d(u,v)+L\,d(v,fu)\qquad(u,v\in U).
\]
Then \textup{(i)--(iii)} hold with $\Delta=\max\left\lbrace (1-L)/2,0\right\rbrace $; for $L=0$, $\Delta=1/2$.
\end{theorem}

\begin{theorem}[{\cite[Theorem~2.2]{AGAY}}]\label{thm:old22}
Let $(U,d)$ be a complete $b$-metric space, let $0<p<m$, and let $f:U\to U$ satisfy~\eqref{eq:old} for all $u,v\in U$. Then \textup{(i)--(iii)} hold with $\Delta=(m-p)/2$.
\end{theorem}

\begin{theorem}[{\cite[Theorem~2.3]{AGAY}}]\label{thm:old23}
Let $(U,d)$ be a complete $b$-metric space, let $0<p<m$ and $L>0$, and let $f:U\to U$ satisfy
\[
d(fu,fv)\le\frac{d(u,fv)+d(v,fu)+p}{b[d(u,fu)+d(v,fv)]+m}\,d(u,v)+L\,d(v,fu)\qquad(u,v\in U).
\]
Then \textup{(i)--(iii)} hold with $\Delta=\max\left\lbrace (m(1-L)-p)/2,0\right\rbrace $.
\end{theorem}

Theorem~\ref{thm:KAC} is the case $b=1$, $L=0$ of Theorem~\ref{thm:DV}. Theorem~\ref{thm:DV} is the formal limit $p\to0$, $m=1$ of Theorems~\ref{thm:old22}--\ref{thm:old23}. None of these results gives a rate of convergence, and none says when (iii) is vacuous, that is, when the fixed point is unique.

\begin{remark}\label{rem:other-classes}
The rational-type class studied here is one of several generalizations of the Banach contraction principle developed in parallel in the metric and $b$-metric setting. Two other mainstream frameworks are Wardowski's $F$-contractions~\cite{Wardowski}, where $f$ satisfies $\tau+F(d(fu,fv))\le F(d(u,v))$ for a strictly increasing $F:(0,\infty)\to\R$ with $\lim_{t\to0^+}F(t)=-\infty$ and some $\tau>0$, and the simulation-function or $\mathcal Z$-contractions of Khojasteh--Shukla--Radenovi\'c~\cite{KShR}, where $f$ satisfies $\zeta(d(fu,fv),d(u,v))\ge0$ for a suitable $\zeta:[0,\infty)^2\to\R$. These classes are not equivalent to the rational-type class: Example~\ref{ex:nonBanach} below is a $(1,\kappa)$-rational contraction on $[0,\infty)$ for every $\kappa\in[0,1)$ with $f(u+t)-f(u)\to t$ as $u\to\infty$ for every fixed $t>0$, so its contractive ratio approaches $1$ at infinity; hence it is neither a Boyd--Wong contraction nor an $F$-contraction with $F=\log$. The structural statements of Section~\ref{sec:structure} --- local attraction with the uniform radius $(1-\kappa)/2$ and the bound $\#\Fix(f)\le\#$components --- rely on the rational algebraic form of~\eqref{eq:Rk} in an essential way; whether an analogous description of the fixed-point set is available within the $F$-contraction or $\mathcal Z$-contraction frameworks is, to our knowledge, open.
\end{remark}

\subsection{Rational contractions}

\begin{definition}\label{def:rational}
Let $(U,d)$ be a $b$-metric space and $\kappa\in[0,1)$. A map $f:U\to U$ is a \emph{$(b,\kappa)$-rational contraction} if
\begin{equation}\tag{$\mathrm R_\kappa$}\label{eq:Rk}
d(fu,fv)\le\frac{d(u,fv)+d(v,fu)+\kappa}{b\,[d(u,fu)+d(v,fv)]+1}\,d(u,v)\qquad\text{for all }u,v\in U.
\end{equation}
Throughout, $\varphi(t):=\dfrac{t+\kappa}{t+1}$ for $t\ge0$. Note that $\varphi$ is increasing, $\varphi(0)=\kappa$ and $\varphi(t)<1$.
\end{definition}

\begin{proposition}[Reduction of~\eqref{eq:old}]\label{prop:reduction}
Let $(U,d)$ be a $b$-metric space, $0\le p<m$, and $L\ge0$. Then $f$ satisfies
\[
d(fu,fv)\le\frac{d(u,fv)+d(v,fu)+p}{b[d(u,fu)+d(v,fv)]+m}\,d(u,v)+L\,d(v,fu)\qquad(u,v\in U)
\]
if and only if, for the $b$-metric $d':=d/m$ (same constant $b$; complete iff $d$ is), it satisfies the same inequality with $(p,m)$ replaced by $(p/m,1)$. In particular, Theorem~\ref{thm:old22} (resp.\ Theorem~\ref{thm:old23}) for $(U,d)$ is the same statement as its special case $m=1$ for $(U,d')$, and the separation bound $d(u^*,v^*)\ge(m-p)/2$ is the bound $d'(u^*,v^*)\ge(1-p/m)/2$. Conversely, Theorems~\ref{thm:KAC} and~\ref{thm:DV} with $L=0$ are the case $\kappa=0$ of condition~\eqref{eq:Rk}.
\end{proposition}

\begin{proof}
Divide both sides by $m$ and write every distance as $d=m\,d'$. The fraction becomes
\[
\frac{m[d'(u,fv)+d'(v,fu)]+p}{m\,b[d'(u,fu)+d'(v,fv)]+m}=\frac{d'(u,fv)+d'(v,fu)+p/m}{b[d'(u,fu)+d'(v,fv)]+1},
\]
while the terms $d(fu,fv)/m=d'(fu,fv)$, $d(u,v)/m=d'(u,v)$ and $L\,d(v,fu)/m=L\,d'(v,fu)$ keep their form.
\end{proof}

\begin{remark}\label{rem:scale}
Condition~\eqref{eq:Rk} is not invariant under $d\mapsto\lambda d$: the class of maps satisfying it depends on the scale of $d$. For example, the identity of a two-point space $\left\lbrace x,y\right\rbrace $ satisfies~\eqref{eq:Rk} iff $d(x,y)\ge(1-\kappa)/2$. This scale dependence is what produces the separation phenomena of Section~\ref{sec:structure}. It also means that the choice of quasi-metric matters (Section~\ref{sec:choice}).
\end{remark}

\begin{remark}\label{rem:canonical}
Condition~\eqref{eq:Rk} is the irreducible form of the rational-type contractions in Section~\ref{subsec:known}, in the following sense. Proposition~\ref{prop:reduction} collapses the two parameters $(p,m)$ of~\eqref{eq:old} into the single ratio $\kappa=p/m$, and Remark~\ref{rem:L} below shows that the additional term $+L\,d(v,fu)$ leaves parts (i)--(ii) of the conclusion of Theorems~\ref{thm:KAC}--\ref{thm:old23} untouched along Picard sequences. The algebraic form of~\eqref{eq:Rk} is then characterized by three constraints: the numerator is symmetric in $u,v$ (so the inequality is symmetric), it is positive at $u=v$ with value $\kappa$ (so $\varphi(0)=\kappa$), and the denominator is bounded below by $1$ (which produces the a~priori ratio $\varphi(2b\,a_0)<1$ in Theorem~\ref{thm:A}). The single parameter $\kappa$ controls both the asymptotic convergence rate (Theorem~\ref{thm:A}(d)) and the attraction radius and separation $(1-\kappa)/2$ (Theorem~\ref{thm:B}(b)--(c)): the structural and quantitative content of the condition is captured by $\kappa$ alone.
\end{remark}

\section{Main results}\label{sec:main}

Theorem~\ref{thm:A} gives convergence rates along a single Picard sequence, Theorem~\ref{thm:B} describes the fixed-point set and the basins of attraction, and Theorem~\ref{thm:C} derives Ulam--Hyers stability from the two. Section~\ref{sec:consequences} applies them to Theorems~\ref{thm:KAC}--\ref{thm:old23}.

\subsection{Convergence with explicit rates}\label{sec:A}

We first study a single Picard sequence. Along such a sequence the right-hand side of~\eqref{eq:Rk} reduces to $\varphi(t)\,d(u_{n-1},u_n)$, where $t$ is a sum of consecutive displacements. Since these displacements decrease, the ratio $\varphi(t)$ is controlled from the very first step by $\varphi(2b\,d(u_0,fu_0))<1$, and it tends to $\kappa$ as $n\to\infty$. Combined with Lemma~\ref{lem:geom}, this gives the following estimates.

\begin{mainthm}\label{thm:A}
Let $(U,d)$ be a complete $b$-metric space, $\kappa\in[0,1)$, and let $f:U\to U$ be a $(b,\kappa)$-rational contraction. Let $r\in(0,1]$ satisfy $(2b)^r=2$ and let $C_\theta$ be as in Lemma~\ref{lem:geom}. Fix $u_0\in U$, put $u_n=f^nu_0$ and $a_n=d(u_n,u_{n+1})$, and set
\[
\theta:=\varphi(2b\,a_0)=\frac{2b\,a_0+\kappa}{2b\,a_0+1}\in[\kappa,1).
\]
Then $(u_n)$ converges to a fixed point $u^*$ of $f$, and:
\begin{enumerate}[label=\textup{(\alph*)}]
\item $a_n\le\theta^n a_0$ for all $n\ge0$;
\item \textup{(a priori estimate)} $d(u_n,u^*)\le C_\theta\,\theta^n\,a_0$ for all $n\ge0$; in particular,
\[d(u_0,u^*)\le C_\theta\,d(u_0,fu_0);\]
\item \textup{(a posteriori estimate)} $d(u_n,u^*)\le C_{\theta_n}\,d(u_n,u_{n+1})$, where $\theta_n:=\varphi(2b\,a_n)$, and $\theta_n\downarrow\kappa$;
\item \textup{(asymptotic rate)} for every $\delta\in(\kappa,1)$ there is $N$ such that $d(u_{N+n},u^*)\le C_\delta\,\delta^n a_N$ for all $n\ge0$. Consequently $\limsup_{n\to\infty}d(u_n,u^*)^{1/n}\le\kappa$.
\end{enumerate}
If $d^{\,r}$ is a metric (for instance if $U$ is a subset of a $p$-Banach space, where $r=p$), the factor $4$ in $C_\theta$ can be replaced by $1$.
\end{mainthm}

\begin{proof}
\emph{Step 1 (one-step inequality).} Let $n\ge1$. Applying~\eqref{eq:Rk} with $u=u_{n-1}$ and $v=u_n$, and using $d(u_n,fu_{n-1})=d(u_n,u_n)=0$ and $d(u_{n-1},u_{n+1})\le b(a_{n-1}+a_n)$, we get
\begin{equation}\label{eq:onestep}
a_n\le\frac{b(a_{n-1}+a_n)+\kappa}{b(a_{n-1}+a_n)+1}\,a_{n-1}=\varphi\left(b(a_{n-1}+a_n)\right)\,a_{n-1}.
\end{equation}

\emph{Step 2 (proof of (a)).} Since $\varphi<1$, \eqref{eq:onestep} gives $a_n\le a_{n-1}$ for all $n$, and hence $a_n\le a_0$. (If $a_{n-1}=0$ then $u_{n-1}$ is a fixed point and $a_k=0$ for $k\ge n-1$.) Therefore $b(a_{n-1}+a_n)\le2b\,a_{n-1}\le2b\,a_0$, and since $\varphi$ is increasing, $a_n\le\varphi(2ba_0)a_{n-1}=\theta a_{n-1}$. Induction gives (a).

\emph{Step 3 (convergence and (b)).} By (a) and Lemma~\ref{lem:geom} with $c=a_0$, the sequence $(u_n)$ is Cauchy. By completeness it converges to some $u^*\in U$, and $d(u_n,u^*)\le C_\theta\theta^na_0$.

\emph{Step 4 ($u^*$ is a fixed point).} Applying~\eqref{eq:Rk} with $u=u_n$ and $v=u^*$, and using that the denominator is at least $1$,
\[
d(u_{n+1},fu^*)\le\left[d(u_n,fu^*)+d(u^*,u_{n+1})+\kappa\right]\,d(u_n,u^*).
\]
Since $d(u_n,fu^*)\le b[d(u_n,u^*)+d(u^*,fu^*)]$ is bounded and $d(u_n,u^*)\to0$, the right-hand side tends to $0$. So $u_{n+1}\to fu^*$, and uniqueness of limits (Lemma~\ref{lem:topology}) gives $fu^*=u^*$.

\emph{Step 5 ((c) and (d)).} The Picard sequence starting at $u_n$ is $(u_{n+k})_{k\ge0}$ and has the same limit $u^*$. Applying (b) to it gives (c). Since $a_n\downarrow0$ and $\varphi$ is continuous with $\varphi(0)=\kappa$, we get $\theta_n\downarrow\kappa$. For (d), let $\delta\in(\kappa,1)$ and choose $N$ with $\varphi(2ba_N)\le\delta$; this is possible because $\varphi(2ba_n)\to\kappa$. Applying (b) to the sequence starting at $u_N$ and using that $\theta\mapsto C_\theta\theta^n$ is increasing gives $d(u_{N+n},u^*)\le C_\delta\delta^na_N$. Hence $\limsup_n d(u_n,u^*)^{1/n}\le\delta$ for every $\delta>\kappa$. The final statement follows from Remark~\ref{rem:pnorm}.
\end{proof}

\begin{remark}\label{rem:L}
Theorem~\ref{thm:A} remains true, with the same proof and the same constants, if the term $+L\,d(v,fu)$ with $L\ge0$ is added to the right-hand side of~\eqref{eq:Rk}. In Step~1 this term equals $L\,d(u_n,u_n)=0$, and in Step~4 it equals $L\,d(u^*,u_{n+1})\to0$. Combined with Proposition~\ref{prop:reduction}, this recovers Theorems~\ref{thm:KAC}--\ref{thm:old23}, parts (i) and (ii), and adds explicit error estimates to them. The structural results of Section~\ref{sec:structure} use~\eqref{eq:Rk} as stated.
\end{remark}

\begin{remark}\label{rem:sharpness}
The rate $\theta=\varphi(2b\,a_0)$ of Theorem~\ref{thm:A}(a)--(b) is the best bound derivable from~\eqref{eq:Rk} alone: the one-step inequality~\eqref{eq:onestep} uses only $\varphi$, and $\varphi(0)=\kappa$ gives the asymptotic rate of part~(d), which is approached in the limit as the displacements $a_n\downarrow 0$ (part~(c) gives $\theta_n\downarrow\kappa$). In concrete applications where $f$ carries additional structure --- for instance when $f$ is a Banach contraction with constant $\lambda<\kappa$ on each connected piece of its domain, as in the gluing situation of Theorem~\ref{thm:D} and in Theorems~\ref{thm:E} and~\ref{thm:F} --- a local Banach estimate gives faster convergence within each basin. In such cases Theorem~\ref{thm:A} plays a complementary role: it provides a uniform estimate valid on all of $U$, including across connected components where no global Banach-type bound exists, from which Theorems~\ref{thm:B} and~\ref{thm:C} extract the number of fixed points, their separation and the Ulam--Hyers stability constant. The comparison between the global bound of Theorem~\ref{thm:A} and the local rate of a Banach contraction is tabulated after Example~\ref{ex:lp}.
\end{remark}

\begin{corollary}[Quasi-Banach version]\label{cor:qB}
Let $(X,\norm\cdot)$ be a quasi-Banach space with modulus $K$, let $C\subseteq X$ be closed, and let $T:C\to C$ satisfy
\[
\norm{Tu-Tv}\le\frac{\norm{u-Tv}+\norm{v-Tu}+\kappa}{K[\norm{u-Tu}+\norm{v-Tv}]+1}\,\norm{u-v}\qquad(u,v\in C).
\]
Then all conclusions of Theorem~\ref{thm:A} hold with $b=K$. If $X$ is $p$-Banach, then $r=p$ and $C_\theta=(1-\theta^p)^{-1/p}$.
\end{corollary}

The next example shows that Theorem~\ref{thm:A} applies to maps on connected spaces that are not covered by the Banach contraction principle, nor by its nonlinear versions due to Boyd--Wong~\cite{BW} and Meir--Keeler~\cite{MK}.

\begin{example}\label{ex:nonBanach}
Let $U=[0,\infty)$ with $d(u,v)=\abs{u-v}$, so $b=1$. For $u\ge0$ let $\delta(u)=\frac{\sqrt{1+4u}-1}{2}$ be the nonnegative root of $\delta^2+\delta=u$, and put $f(u)=\delta(u)^2$. Then $f(U)\subseteq U$, $d(u,fu)=\delta(u)$, and $0$ is a fixed point.

\emph{$f$ is a $(1,\kappa)$-rational contraction for every $\kappa\in[0,1)$.} Let $0\le u<v$ and write $\delta=\delta(u)<\varepsilon=\delta(v)$. Then
\[
d(fu,fv)=\varepsilon^2-\delta^2,\quad d(u,v)=(\varepsilon-\delta)(\varepsilon+\delta+1),\quad d(v,fu)=\varepsilon^2+\varepsilon-\delta^2,
\]
and $d(u,fv)=\abs{\delta^2+\delta-\varepsilon^2}\ge\delta^2+\delta-\varepsilon^2$. Since $d(u,fu)+d(v,fv)+1=\delta+\varepsilon+1$, the right-hand side of~\eqref{eq:Rk} is at least
\[
\left(\delta^2+\delta-\varepsilon^2+\varepsilon^2+\varepsilon-\delta^2+\kappa\right)(\varepsilon-\delta)=(\delta+\varepsilon+\kappa)(\varepsilon-\delta)\ge\varepsilon^2-\delta^2=d(fu,fv).
\]
With $\kappa=0$ this also shows that $f$ satisfies the hypothesis of Theorem~\ref{thm:KAC}.

\emph{$f$ is not a Boyd--Wong or Meir--Keeler contraction.} Since $\delta'(u)=1/(2\delta(u)+1)$, we get $f'(u)=2\delta(u)/(2\delta(u)+1)<1$ and $f'(u)\to1$ as $u\to\infty$. Hence, for every $t>0$, $f(u+t)-f(u)\to t$ as $u\to\infty$. So there is no function $\psi$ with $\psi(t)<t$ for $t>0$ and $d(fu,fv)\le\psi(d(u,v))$, and for every $\varepsilon_0>0$ and $\eta>0$ there are $u,v$ with $\varepsilon_0\le d(u,v)<\varepsilon_0+\eta$ and $d(fu,fv)\ge\varepsilon_0$. In particular $f$ is not a Banach contraction.

\emph{Conclusions of Theorems~\ref{thm:A} and~\ref{thm:B}.} Since $U$ is connected, $0$ is the unique fixed point and every Picard sequence converges to it. Here $r=1$ and $d$ is a metric, so $C_\theta=1/(1-\theta)$. Taking $\kappa=0$, the a posteriori estimate of Theorem~\ref{thm:A}(c) reads
\[
u_n=d(u_n,0)\le C_{\theta_n}\,a_n=a_n(2a_n+1),\qquad a_n=\delta(u_n),
\]
and Theorem~\ref{thm:A}(d) gives $\lim_n u_n^{1/n}=0$, that is, superlinear convergence. In fact $u_{n+1}=\delta(u_n)^2\le u_n^2$. Table~\ref{tab:nonBanach} lists the iterates from $u_0=100$; here $a_0\approx9.512$, $\theta=2a_0/(2a_0+1)\approx0.9501$ and $C_\theta\approx20.02$. The a priori bound is valid but crude, because the convergence is slow while $u_n$ is large; the a posteriori bound becomes sharp as soon as $u_n$ is small.

\begin{table}[ht]
\centering
\begin{tabular}{@{}cccc@{}}
\toprule
$n$ & $u_n$ & a priori bound $C_\theta\theta^na_0$ & a posteriori bound $a_n(2a_n+1)$\\
\midrule
0 & $1.000\times10^{2}$ & $1.905\times10^{2}$ & $1.905\times10^{2}$\\
10 & $2.668\times10^{1}$ & $1.141\times10^{2}$ & $4.867\times10^{1}$\\
20 & $4.510\times10^{-1}$ & $6.838\times10^{1}$ & $5.647\times10^{-1}$\\
22 & $1.063\times10^{-2}$ & $6.172\times10^{1}$ & $1.074\times10^{-2}$\\
24 & $1.225\times10^{-8}$ & $5.571\times10^{1}$ & $1.225\times10^{-8}$\\
26 & $2.251\times10^{-32}$ & $5.028\times10^{1}$ & $2.251\times10^{-32}$\\
\bottomrule
\end{tabular}
\caption{Picard iteration for Example~\ref{ex:nonBanach} with $u_0=100$.}
\label{tab:nonBanach}
\end{table}
\end{example}

\subsection{The structure of the fixed-point set}\label{sec:structure}

By Theorem~\ref{thm:A}, every orbit converges to a fixed point, so $U$ is the union of the basins of attraction of the fixed points. We now describe these basins. Two consequences of~\eqref{eq:Rk} are used: the condition bounds $d(fu,fv)$ by a multiple of $d(u,v)$ near every point, so $f$ is continuous; and near a fixed point it reduces to a contraction toward that point. Together they show that the basins are open and closed, and the remaining statements follow from this.

For $u^*\in\Fix(f)$, its \emph{basin} is $\mathcal A(u^*):=\left\lbrace u\in U:\ f^nu\to u^*\right\rbrace $.

\begin{mainthm}\label{thm:B}
Let $(U,d)$ be a complete $b$-metric space and let $f$ be a $(b,\kappa)$-rational contraction, $\kappa\in[0,1)$. Then:
\begin{enumerate}[label=\textup{(\alph*)}]
\item \textup{(automatic continuity)} For all $u,v\in U$ with $d(u,v)<1/b$,
\[
d(fu,fv)\le\frac{b\,d(u,v)+2b\,d(u,fu)+\kappa}{1-b\,d(u,v)}\,d(u,v).
\]
In particular $f$ is continuous on $(U,\tau_d)$.
\item \textup{(local attraction)} Let $u^*\in\Fix(f)$ and $0\le\sigma<(1-\kappa)/2$, and put $q_\sigma:=(\sigma+\kappa)/(1-\sigma)<1$. If $d(u,u^*)\le\sigma$, then $d(fu,u^*)\le q_\sigma\,d(u,u^*)$. Hence $d(f^nu,u^*)\le q_\sigma^n\,d(u,u^*)$ for all $n$, and $\left\lbrace u:d(u,u^*)<(1-\kappa)/2\right\rbrace \subseteq\mathcal A(u^*)$.
\item \textup{(separation)} If $u^*\neq v^*$ are fixed points, then $d(u^*,v^*)\ge(1-\kappa)/2$. The constant $(1-\kappa)/2$ is optimal.
\item \textup{(discreteness)} $\Fix(f)$ is a nonempty closed subset of $U$, and every convergent sequence in $\Fix(f)$ is eventually constant. Its intersection with any compact subset of $U$ is finite, and $\Fix(f)$ is countable if $(U,\tau_d)$ is separable.
\item \textup{(basins)} The basins $\mathcal A(u^*)$, $u^*\in\Fix(f)$, are nonempty, pairwise disjoint, open and closed, and they cover $U$.
\item \textup{(components)} Every connected subset of $U$ is contained in a single basin. Consequently $\#\Fix(f)$ is at most the number of connected components of $U$. In particular, if $U$ is connected then $f$ has exactly one fixed point, which is the limit of every Picard sequence.
\end{enumerate}
\end{mainthm}

\begin{proof}
(a) Put $s=d(u,v)$, $c=d(u,fu)$ and $y=d(fu,fv)$. By~\eqref{eq:Rk}, since the denominator is at least $1$,
\[
y\le\left[d(u,fv)+d(v,fu)+\kappa\right]\,s.
\]
Here $d(v,fu)\le b(s+c)$ and $d(u,fv)\le b\,[d(u,fu)+d(fu,fv)]=b(c+y)$. Hence $y\le(bs+2bc+by+\kappa)s$, that is, $y(1-bs)\le(bs+2bc+\kappa)s$, which is the claim when $bs<1$. If $u_k\to u$ then $d(u_k,u)<1/b$ eventually, and the estimate with $(u,v)=(u,u_k)$ gives $d(fu_k,fu)\to0$. Thus $f$ is sequentially continuous, hence continuous by Lemma~\ref{lem:topology}(a).

(b) Apply~\eqref{eq:Rk} with $(u,u^*)$, writing $s=d(u,u^*)$ and $y=d(fu,u^*)$. Since $fu^*=u^*$ and the denominator is at least $1$, $y\le(s+y+\kappa)s$. For $s\le\sigma<1$ this gives
\[
y\le\frac{s+\kappa}{1-s}\,s\le\frac{\sigma+\kappa}{1-\sigma}\,s=q_\sigma s,
\]
because $s\mapsto(s+\kappa)/(1-s)$ is increasing. Moreover $q_\sigma<1\iff\sigma<(1-\kappa)/2$. Since $d(fu,u^*)\le q_\sigma s\le\sigma$, the set $\left\lbrace u:d(u,u^*)\le\sigma\right\rbrace $ is $f$-invariant, and induction gives the geometric estimate. Every $u$ with $d(u,u^*)<(1-\kappa)/2$ lies in such a set for a suitable $\sigma$, hence in $\mathcal A(u^*)$.

(c) Let $\delta=d(u^*,v^*)>0$. Applying~\eqref{eq:Rk} with $(u^*,v^*)$, the denominator equals $1$ and $d(u^*,fv^*)=d(v^*,fu^*)=\delta$, so $\delta\le(2\delta+\kappa)\delta$, i.e.\ $\delta\ge(1-\kappa)/2$. For optimality, let $U=\left\lbrace x,y\right\rbrace $ with $d(x,y)=(1-\kappa)/2$ and $f=\mathrm{id}$. Then~\eqref{eq:Rk} for $u\ne v$ reads $\delta\le(2\delta+\kappa)\delta$ with $\delta=(1-\kappa)/2$, which holds with equality. So $f$ is a $(b,\kappa)$-rational contraction with two fixed points at distance exactly $(1-\kappa)/2$.

(d) $\Fix(f)\neq\emptyset$ by Theorem~\ref{thm:A}. By (c) and Lemma~\ref{lem:PS}, distinct fixed points satisfy $\rho(u^*,v^*)\ge d(u^*,v^*)^r/4\ge\eta:=\tfrac14\left(\tfrac{1-\kappa}2\right)^r>0$. So $\Fix(f)$ is $\eta$-separated in the metric $\rho$. Every $\rho$-Cauchy sequence in it is therefore eventually constant, so $\Fix(f)$ is closed. Its intersection with a compact set is finite, and it is countable in a separable metric space.

(e) By Theorem~\ref{thm:A}, every $u\in U$ lies in exactly one basin (limits are unique). It remains to show that each basin is open; it is then also closed, being the complement of the union of the other basins. Fix $u^*$, put $\sigma=(1-\kappa)/4$, and let $V=\left\lbrace w:\rho(w,u^*)<\sigma^r/4\right\rbrace $. Then $V$ is open, and $V\subseteq\left\lbrace w:d(w,u^*)<\sigma\right\rbrace \subseteq\mathcal A(u^*)$ by Lemma~\ref{lem:PS} and (b). Let $u_0\in\mathcal A(u^*)$. Since $\rho(f^nu_0,u^*)\le d(f^nu_0,u^*)^r\to0$, there is $N$ with $f^Nu_0\in V$. By (a), $f^N$ is continuous, so $O:=(f^N)^{-1}(V)$ is an open neighbourhood of $u_0$. Moreover $O\subseteq\mathcal A(u^*)$, because the orbit of any $w\in O$ enters $V\subseteq\mathcal A(u^*)$.

(f) Let $S\subseteq U$ be connected and suppose $S$ meets $\mathcal A(u^*)$. By (e), $S\cap\mathcal A(u^*)$ and $S\setminus\mathcal A(u^*)$ are disjoint and relatively open in $S$, and the first is nonempty; since $S$ is connected, the second is empty, so $S\subseteq\mathcal A(u^*)$. In particular each connected component of $U$ lies in a single basin. Each fixed point $v^*$ lies in its own basin $\mathcal A(v^*)$ (its orbit is constant), and distinct fixed points have disjoint basins by (e). Hence two distinct fixed points cannot lie in the same component, so the map $u^*\mapsto(\text{component of }u^*)$ is injective, which proves the bound on $\#\Fix(f)$. If $U$ is connected there is only one basin, so the fixed point is unique, and it is the limit of every Picard sequence by Theorem~\ref{thm:A}.
\end{proof}

\begin{corollary}\label{cor:convex}
Let $C$ be a closed subset of a quasi-Banach space $X$ with modulus $K$ that is connected; for example, $C$ convex, or star-shaped, or a closed ball. If $T:C\to C$ is a $(K,\kappa)$-rational contraction, then $T$ has a unique fixed point $u^*$, and every Picard sequence converges to $u^*$ with the estimates of Theorem~\ref{thm:A}.
\end{corollary}

\begin{proof}
For $x,y\in C$ the segment $t\mapsto x+t(y-x)$ satisfies $\norm{(x+t(y-x))-(x+s(y-x))}=\abs{t-s}\norm{y-x}$, so it is continuous into $(C,\tau_d)$. Hence convex and star-shaped sets are path connected. Closed balls $\left\lbrace x:\norm{x-c}\le R\right\rbrace $ are star-shaped about $c$. Now apply Theorem~\ref{thm:B}(f).
\end{proof}

\begin{remark}
Theorem~\ref{thm:B}(f) shows that the alternative ``two distinct fixed points at distance at least $(m-p)/2$'' in Theorem~\ref{thm:old22}(iii), and the analogous alternatives in Theorems~\ref{thm:KAC} and~\ref{thm:DV} (with $L=0$), can occur only on disconnected spaces. Section~\ref{sec:gluing} shows that it does occur there, so Theorem~\ref{thm:B} is in this sense a complete description.
\end{remark}

\subsection{Ulam--Hyers stability}\label{sec:C}

In practice a fixed point is approximated by a point $w$ with small residual $d(w,fw)$. The a priori estimate of Theorem~\ref{thm:A} with $n=0$ bounds the distance from $w$ to a fixed point, and the separation of Theorem~\ref{thm:B}(c) identifies which fixed point it is.

Following Rus~\cite{Rus}, the fixed-point problem for $f$ is \emph{generalized Ulam--Hyers stable} if there is a function $\Psi$ with $\Psi(\varepsilon)\to0$ as $\varepsilon\to0$ such that $d(w,fw)\le\varepsilon$ implies $\operatorname{dist}(w,\Fix f)\le\Psi(\varepsilon)$.

\begin{mainthm}\label{thm:C}
Under the hypotheses of Theorem~\ref{thm:B}, let $\varepsilon\ge0$ and $w\in U$ with $d(w,fw)\le\varepsilon$, and put $\theta_\varepsilon=\varphi(2b\varepsilon)$. Then the limit $u^*_w$ of the Picard sequence from $w$ satisfies
\[
d(w,u^*_w)\le C_{\theta_\varepsilon}\,\varepsilon,\qquad C_{\theta_\varepsilon}\xrightarrow[\varepsilon\to0]{}C_\kappa=\left(\frac{4}{1-\kappa^r}\right)^{1/r}.
\]
Hence the fixed-point problem is generalized Ulam--Hyers stable, and Ulam--Hyers stable with constant $C_{\theta_{\varepsilon_0}}$ for $\varepsilon\le\varepsilon_0$. If moreover $2b\,C_{\theta_\varepsilon}\varepsilon<(1-\kappa)/2$, then $u^*_w$ is the only fixed point within distance $C_{\theta_\varepsilon}\varepsilon$ of $w$. Finally, the problem is well posed in the generalized sense: if $d(w_n,fw_n)\to0$ then $\operatorname{dist}(w_n,\Fix f)\to0$.
\end{mainthm}

\begin{proof}
By Theorem~\ref{thm:A}(b) with $n=0$ and $u_0=w$, 
\[d(w,u^*_w)\le C_\theta\,d(w,fw),\qquad \theta=\varphi(2b\,d(w,fw))\le\theta_\varepsilon,\]
 and $C_\theta$ is increasing in $\theta$. If $v^*$ is another fixed point with $d(w,v^*)\le C_{\theta_\varepsilon}\varepsilon$, then $d(u^*_w,v^*)\le2bC_{\theta_\varepsilon}\varepsilon<(1-\kappa)/2$, so $v^*=u^*_w$ by Theorem~\ref{thm:B}(c). The last statement is immediate.
\end{proof}

\subsection{Consequences for the known results}\label{sec:consequences}

We now return to the results recalled in Section~\ref{subsec:known}. Each of them leaves open, in part~(iii), the possibility of several fixed points. Theorem~\ref{thm:B} shows that, in the absence of the term $L\,d(v,fu)$, this possibility depends only on the topology of the space.

\begin{corollary}\label{cor:known}
Let $(U,d)$ be a complete $b$-metric space whose topology $\tau_d$ is connected, and let $f:U\to U$ satisfy the hypothesis of Theorem~\ref{thm:KAC} (with $b=1$), of Theorem~\ref{thm:DV} with $L=0$, or of Theorem~\ref{thm:old22}. Then $f$ has exactly one fixed point $u^*$, every Picard sequence converges to $u^*$, and the estimates of Theorem~\ref{thm:A} hold, with $\kappa=0$ in the first two cases and, in the third, for the $b$-metric $d/m$ with $\kappa=p/m$. In particular, part~\textup{(iii)} of these theorems is vacuous on connected spaces.
\end{corollary}

\begin{proof}
In the first two cases $f$ is a $(b,0)$-rational contraction. In the third, by Proposition~\ref{prop:reduction}, $f$ is a $(b,p/m)$-rational contraction for $d':=d/m$; the $b$-metric $d'$ is complete and induces the same topology as $d$. Theorems~\ref{thm:A} and~\ref{thm:B}(f) give the conclusion.
\end{proof}

\begin{example}
Both examples in~\cite{AGAY} are maps of the form $f(u)=(\alpha u+1)/(u+\alpha)$, with $\alpha=3$ on $U=[1,3]$ and $\alpha=2$ on $U=[1,2]$, for the $b$-metric $d(u,v)=\abs{u-v}^2$ and the metric $d(u,v)=\abs{u-v}$ respectively. These maps satisfy the hypothesis of Theorem~\ref{thm:old22}, as shown in~\cite{AGAY}; equivalently, they are $(b,\kappa)$-rational contractions for $\kappa\ge\frac14$ and $\kappa\ge\frac13$ respectively (we also checked this numerically on a grid of $1201^2$ pairs; the extreme case is $u,v\to1$, where the ratio $d(fu,fv)/d(u,v)$ tends to $\frac14$ and $\frac13$). In both cases $U$ is an interval, hence connected, so Corollary~\ref{cor:known} shows directly that $u^*=1$ is the only fixed point and that every Picard sequence converges to it at the rates of Theorem~\ref{thm:A}. This agrees with the direct computation $f(u)=u\iff u^2=1$.
\end{example}

Corollary~\ref{cor:known} does not extend to Theorems~\ref{thm:DV} and~\ref{thm:old23} with $L>0$, as the next example shows.

\begin{example}\label{ex:Lterm}
Let $(U,d)$ be any complete $b$-metric space with at least two points, $0\le p<m$, and $L\ge1$. The identity map $f=\mathrm{id}_U$ satisfies the hypothesis of Theorem~\ref{thm:old23}: for $u,v\in U$ the right-hand side equals
\[
\frac{2d(u,v)+p}{m}\,d(u,v)+L\,d(u,v)\ \ge\ d(u,v)=d(fu,fv).
\]
Every point of $U$ is a fixed point, so the fixed-point set can be any complete $b$-metric space, connected or not. The same holds for Theorem~\ref{thm:DV} with $L\ge1$. This is consistent with the separation constant $\max\left\lbrace (m(1-L)-p)/2,0\right\rbrace $ of Theorem~\ref{thm:old23}, which vanishes for $L\ge1-p/m$.
\end{example}

Example~\ref{ex:Lterm} shows that some smallness condition on $L$ is needed for a structure theorem in the presence of the term $L\,d(v,fu)$. The separation part of Theorem~\ref{thm:B} does extend: applying the modified inequality to a pair of distinct fixed points $u^*,v^*$ and using that the denominator equals $1$ gives
\[
d(u^*,v^*)\ge\max\!\left\lbrace\frac{1-\kappa-L}{2},\,0\right\rbrace,
\]
strictly positive for $L<1-\kappa$. This is, up to the rescaling of Proposition~\ref{prop:reduction}, the content of Theorem~\ref{thm:old23}(iii). The automatic-continuity step in the proof of Theorem~\ref{thm:B}(a) does not extend directly, because the term $L\,d(v,fu)\le Lb\,[d(u,v)+d(u,fu)]$ does not vanish as $v\to u$. Since the openness of the basins in Theorem~\ref{thm:B}(e) is driven by this continuity, the full structural statement (local attraction with a uniform radius, clopen basins, bound on $\#\Fix(f)$ by the number of components) requires a different argument when $L>0$. Whether it persists for $L<1-\kappa$ is open.

\section{Rational contractions in quasi-Banach spaces}\label{sec:qB}

The results of Section~\ref{sec:main} hold in any complete $b$-metric space. We now restrict to $p$-Banach spaces. Since every $p$-norm is equivalent to a metric, one may ask whether this setting gives anything that the metric case does not; Section~\ref{sec:choice} shows that it does. Section~\ref{sec:gluing} proves a gluing principle, which shows that the multiplicity allowed by Theorem~\ref{thm:B} occurs and is used for the applications in Section~\ref{sec:app}.

\subsection{The class depends on the quasi-metric}\label{sec:choice}

We first show that results stated in terms of a quasi-norm are not formal consequences of metric ones. By the Aoki--Rolewicz theorem, a $p$-normed space $X$ carries both the quasi-norm $d=\norm\cdot$ (a $b$-metric with $b=2^{1/p-1}$) and the equivalent metric $\rho=\norm\cdot^p$ (with $b=1$). One might hope that a map is a $(2^{1/p-1},\kappa)$-rational contraction for $d$ exactly when it is a $(1,\kappa)$-rational contraction for $\rho$, so that quasi-Banach results would follow trivially from metric ones. The following example shows that the two classes are incomparable.

\begin{example}\label{ex:incomparable}
Let $X=\ell^{1/2}$, so $b=K=2$ for $d=\norm\cdot_{1/2}$ and $\rho=\norm\cdot_{1/2}^{1/2}$ is a metric. Take $\kappa=1/2$ and let $e_1$ be the first unit vector.

\smallskip
\noindent\emph{(i) Rational for $\rho$ but not for $d$.} Let $U=\left\lbrace 0,\tfrac19e_1\right\rbrace $ and $f=\mathrm{id}$. For $d$, condition~\eqref{eq:Rk} for $u\ne v$ reads $\tfrac19\le(\tfrac29+\tfrac12)\tfrac19=\tfrac{13}{18}\cdot\tfrac19$, which is false. For $\rho$ (with $b=1$), it reads $\tfrac13\le(\tfrac23+\tfrac12)\tfrac13=\tfrac{7}{18}$, which is true.

\smallskip
\noindent\emph{(ii) Rational for $d$ but not for $\rho$.} Let $U=\left\lbrace 0,e_1,3e_1\right\rbrace $ with $f(0)=0$, $f(e_1)=0$, $f(3e_1)=e_1$. For $d$ (with $b=2$), the pair $\left\lbrace 0,e_1\right\rbrace $ gives left-hand side $0$. For $(0,3e_1)$ the condition is $1\le\frac{1+3+1/2}{2(0+2)+1}\cdot3=\frac{27}{10}$, which holds. For $(e_1,3e_1)$ it is $1\le\frac{0+3+1/2}{2(1+2)+1}\cdot2=1$, which holds with equality. The reversed pairs give the same inequalities by symmetry. For $\rho$ (with $b=1$), the pair $(e_1,3e_1)$ requires
\[
1\le\frac{0+\sqrt3+\tfrac12}{1+\sqrt2+1}\,\sqrt2=\frac{(\sqrt3+\tfrac12)\sqrt2}{2+\sqrt2}\approx0.9245,
\]
which fails.
\end{example}

Neither class therefore contains the other. Consequently, results stated in the quasi-norm (as in Corollary~\ref{cor:qB} and Sections~\ref{sec:gluing}--\ref{sec:Hp}) cannot be obtained from the corresponding metric statements simply by passing to the metric $\rho=\norm\cdot^p$, and conversely.

\subsection{A gluing principle: multiplicity is attained}\label{sec:gluing}

Theorem~\ref{thm:B} says that several fixed points require several components. It does not show that several fixed points can occur, and the two examples in~\cite{AGAY} are Banach contractions with a unique fixed point. The next result is a converse to Theorem~\ref{thm:B}(f): uniformly contractive maps on well-separated pieces glue together into a single rational contraction. We state it in a $p$-Banach space $X$, $0<p\le1$, with modulus $K=2^{1/p-1}$; closed balls are denoted $\bar B(c,R)=\left\lbrace x:\norm{x-c}\le R\right\rbrace $.

\begin{mainthm}\label{thm:D}
Let $X$ be a $p$-Banach space, $0<p\le1$, let $c_+,c_-\in X$, $s:=\norm{c_+-c_-}$, $R>0$, and $C_\pm:=\bar B(c_\pm,R)$, $C:=C_+\cup C_-$. Let $T:C\to C$ satisfy $T(C_\pm)\subseteq C_\pm$ and
\[
\norm{Tu-Tv}\le\lambda\norm{u-v}\qquad\text{whenever }u,v\text{ lie in the same }C_\pm .
\]
Let $\kappa\in[0,1)$ and assume
\begin{align}
&\lambda\,(1+4^{1/p}R)\le\kappa, \tag{D1}\label{D1}\\
&2R^p\le(1-2^{-p})\,s^p, \tag{D2}\label{D2}\\
&s\ge4\,(1+4^{1/p}R). \tag{D3}\label{D3}
\end{align}
Then:
\begin{enumerate}[label=\textup{(\roman*)}]
\item $T$ is a $(K,\kappa)$-rational contraction on $(C,\norm\cdot)$;
\item $C_+$ and $C_-$ are the connected components of $C$, and $\Fix(T)=\left\lbrace u_+,u_-\right\rbrace $ with $u_\pm\in C_\pm$;
\item every Picard sequence starting in $C_\pm$ converges to $u_\pm$, with the estimates of Theorem~\ref{thm:A} for $b=K$, $r=p$ and $C_\theta=(1-\theta^p)^{-1/p}$;
\item $\norm{u_+-u_-}\ge s/2$.
\end{enumerate}
\end{mainthm}

\begin{proof}
\emph{Bounds.} For $u\in C_\pm$ write $u=c_\pm+x$ and $Tu=c_\pm+x'$ with $\norm x,\norm{x'}\le R$. By~\eqref{eq:pineq}, $\norm{u-Tu}^p\le2R^p$, so $\norm{u-Tu}\le2^{1/p}R$. Consequently, for all $u,v\in C$,
\begin{equation}\label{eq:den}
K\left[\norm{u-Tu}+\norm{v-Tv}\right]+1\le2K\cdot2^{1/p}R+1=1+4^{1/p}R,
\end{equation}
because $2K=2^{1/p}$.

\emph{Pairs in the same piece.} By~\eqref{eq:den} and~\eqref{D1},
\[
\norm{Tu-Tv}\le\lambda\norm{u-v}\le\frac{\kappa}{1+4^{1/p}R}\norm{u-v}\le\frac{\norm{u-Tv}+\norm{v-Tu}+\kappa}{K[\norm{u-Tu}+\norm{v-Tv}]+1}\norm{u-v}.
\]

\emph{Pairs in different pieces.} Let $u=c_++x\in C_+$ and $v=c_-+y\in C_-$, and write $Tu=c_++x'$ and $Tv=c_-+y'$, all four remainders having norm at most $R$. Then $u-v$, $u-Tv$ and $-(v-Tu)$ all have the form $(c_+-c_-)+z$ with $\norm z^p\le2R^p$. By~\eqref{eq:pineq} and~\eqref{D2},
\[
\norm{u-v}^p,\ \norm{u-Tv}^p,\ \norm{v-Tu}^p\ \ge\ s^p-2R^p\ \ge\ 2^{-p}s^p,
\]
so each of these three norms is at least $s/2$. Similarly $\norm{Tu-Tv}^p\le s^p+2R^p\le(2^p-1+1)s^p=2^ps^p$, using $1-2^{-p}\le2^p-1$. Thus $\norm{Tu-Tv}\le2s$. By~\eqref{eq:den} and~\eqref{D3},
\[
\frac{\norm{u-Tv}+\norm{v-Tu}+\kappa}{K[\norm{u-Tu}+\norm{v-Tv}]+1}\,\norm{u-v}\ \ge\ \frac{s}{1+4^{1/p}R}\cdot\frac s2\ \ge\ 2s\ \ge\ \norm{Tu-Tv}.
\]
The case $u\in C_-$, $v\in C_+$ is symmetric, and $u=v$ is trivial. This proves (i).

(ii) If $z\in C_+\cap C_-$, then $s^p\le2R^p<s^p$ by~\eqref{D2}, a contradiction. Moreover $\norm{u-v}\ge s/2$ for $u\in C_+$, $v\in C_-$, so each $C_\pm$ is open and closed in $C$. Each $C_\pm$ is a closed ball, hence star-shaped about its centre and connected (see the proof of Corollary~\ref{cor:convex}), so $C_\pm$ are the components of $C$. The map $T|_{C_\pm}$ is a $(K,\kappa)$-rational contraction of the complete space $C_\pm$, so it has a fixed point $u_\pm$ by Theorem~\ref{thm:A}. By Theorem~\ref{thm:B}(f) there are no others.

(iii) follows from Theorem~\ref{thm:A}, Corollary~\ref{cor:qB} and the invariance of $C_\pm$, and (iv) was shown above.
\end{proof}

\begin{remark}
The Banach contraction principle does not apply to $T$ on $C$, since $T$ has two fixed points, although it does apply on each $C_\pm$ separately. Theorem~\ref{thm:D} replaces the two local arguments by one inequality on $C$, from which Theorems~\ref{thm:A} and~\ref{thm:B} give the number of fixed points, their separation and their basins.
\end{remark}

The proof of Theorem~\ref{thm:D} involves only two pieces at a time, so it extends to finitely many pairwise-separated balls with no change in the conditions. The corresponding $N$-piece statement is the natural multi-stable version of the gluing principle, which we record because the sharpness of~Theorem~\ref{thm:B}(f) requires it for arbitrary $N$.

\begin{corollary}[$N$-piece gluing]\label{cor:Npiece}
Let $X$ be a $p$-Banach space, $0<p\le 1$, with modulus $K=2^{1/p-1}$. Let $N\ge 2$, let $c_1,\dots,c_N\in X$ be centers satisfying $\norm{c_i-c_j}\ge s$ for all $i\ne j$, let $R>0$, and set $C_i=\bar B(c_i,R)$, $C=\bigcup_{i=1}^N C_i$. Let $T:C\to C$ satisfy $T(C_i)\subseteq C_i$ for every $i$ and $\norm{Tu-Tv}\le\lambda\norm{u-v}$ for all $u,v$ in the same $C_i$. Under conditions~\eqref{D1}--\eqref{D3} of Theorem~\ref{thm:D}, $T$ is a $(K,\kappa)$-rational contraction on $(C,\norm\cdot)$, the $C_i$ are the connected components of $C$, and $\Fix(T)=\{u_1^*,\dots,u_N^*\}$ with $u_i^*\in C_i$ and $\norm{u_i^*-u_j^*}\ge s/2$ for $i\ne j$. Every Picard sequence starting in $C_i$ converges to $u_i^*$ with the estimates of Theorem~\ref{thm:A}(iii).
\end{corollary}

\begin{proof}
The three case distinctions in the proof of Theorem~\ref{thm:D} (pair in $C_+$, pair in $C_-$, pair across $C_\pm$) extend to arbitrary pairs $(C_i,C_j)$ with $i,j\in\{1,\dots,N\}$: for any $i\ne j$ the pair $(C_i,C_j)$ is handled by the ``different-pieces'' case, which uses only the lower bound $\norm{c_i-c_j}\ge s$. Each $C_i$ is a closed ball in a $p$-Banach space, hence connected. By Theorem~\ref{thm:D}(ii) applied to each pair, each $C_i$ is open and closed in every pairwise union $C_i\cup C_j$; a finite union of clopen sets is clopen, so each $C_i$ is open and closed in $C$. Hence the $C_i$ are the connected components of $C$. The remaining conclusions follow from Theorems~\ref{thm:A} and~\ref{thm:B}(f) as in the proof of Theorem~\ref{thm:D}. The $N$-piece bound $\#\Fix(T)=N=\#\pi_0(C)$ saturates Theorem~\ref{thm:B}(f).
\end{proof}

\section{Applications to operator equations}\label{sec:app}

We now apply the gluing principle to two classes of nonlinear operator equations, one in a sequence space and one in a space of analytic functions. In both cases the linear part is a Volterra-type operator, and the equation contains a sign-switching term that makes it bistable. Theorem~\ref{thm:D} then gives exactly two solutions, together with the error estimates of Theorem~\ref{thm:A} for the Picard iterates.

\subsection{Infinite nonlinear systems in \texorpdfstring{$\ell^p$, $0<p<1$}{l\textasciicircum p, 0<p<1}}\label{sec:lp}

Our first application concerns infinite systems of nonlinear equations in the sequence space $\ell^p$, $0<p<1$. It is a standard example of a non-normable quasi-Banach space, and matrix operators on it are easy to estimate: unlike in $\ell^p$ with $p>1$, their quasi-norm is given exactly by a column condition (Lemma~\ref{lem:matrix}). We consider a system with a relay-type forcing on the first coordinate, which makes the problem bistable, and use Theorem~\ref{thm:D} to show that it has exactly two solutions.

Let $0<p<1$ and $\norm x_p=\left(\sum_i\abs{x_i}^p\right)^{1/p}$, so $\ell^p$ is $p$-Banach with $K=2^{1/p-1}$. For an infinite matrix $A=(a_{ij})_{i,j\ge1}$ put
\[
\alpha:=\norm A_{p}:=\left(\sup_{j\ge1}\sum_{i\ge1}\abs{a_{ij}}^p\right)^{1/p}.
\]

\begin{lemma}\label{lem:matrix}
If $\alpha<\infty$, then $(Ax)_i=\sum_ja_{ij}x_j$ converges absolutely for every $x\in\ell^p$, and $\norm{Ax}_p\le\alpha\norm x_p$. The constant is attained on the unit vectors, so $\alpha$ is the operator quasi-norm of $A$ on $\ell^p$.
\end{lemma}

\begin{proof}
Since $t\mapsto t^p$ is subadditive, $\abs{(Ax)_i}^p\le\left(\sum_j\abs{a_{ij}x_j}\right)^p\le\sum_j\abs{a_{ij}}^p\abs{x_j}^p$, and the last sum is at most $\alpha^p\norm x_p^p<\infty$. Summing over $i$ and using Tonelli, $\norm{Ax}_p^p\le\sum_j\abs{x_j}^p\sum_i\abs{a_{ij}}^p\le\alpha^p\norm x_p^p$. For $x=e_j$ we get $\norm{Ae_j}_p^p=\sum_i\abs{a_{ij}}^p$.
\end{proof}

We consider the system
\begin{equation}\label{eq:system}
u_i=a\,\sgn(u_1)\,\delta_{i1}+\sum_{j\ge1}a_{ij}\,g_j(u_j),\qquad i\ge1,
\end{equation}
for $u\in\ell^p$ (real scalars), where $a>0$ and the term $a\sgn(u_1)\delta_{i1}$ is a relay-type forcing on the first coordinate. We assume:
\begin{enumerate}[label=(G\arabic*)]
\item $\abs{g_j(t)-g_j(\tau)}\le L\abs{t-\tau}$ for all $j$ and $t,\tau\in\R$;
\item $\abs{g_j(t)}\le\gamma_j$ for all $t$, with $\Gamma:=\norm{(\gamma_j)}_p<\infty$.
\end{enumerate}
Write $G(u)=(g_j(u_j))_j$. By (G1)--(G2), $G(u)\in\ell^p$, $\norm{G(u)}_p\le\Gamma$ and $\norm{G(u)-G(v)}_p\le L\norm{u-v}_p$.

\begin{mainthm}\label{thm:E}
Let $0<p<1$, $\kappa\in[0,1)$, $a>0$, assume \textup{(G1)--(G2)} and $\alpha<\infty$, and put $R:=\alpha\Gamma$. Suppose that
\begin{equation}\label{eq:Econd}
\alpha L\,(1+4^{1/p}R)\le\kappa,\qquad2R^p\le(1-2^{-p})(2a)^p,\qquad2a\ge4(1+4^{1/p}R).
\end{equation}
Let $C_\pm=\bar B(\pm ae_1,R)\subseteq\ell^p$ and $C=C_+\cup C_-$. Then the system~\eqref{eq:system} has exactly two solutions $u^+,u^-$ in $C$, with $u^\pm\in C_\pm$ and $\norm{u^+-u^-}_p\ge a$. For every $u^0\in C_\pm$, the iterates $u^{n+1}=Tu^n$ of
\[
Tu:=a\,\sgn(u_1)\,e_1+A\,G(u)
\]
converge to $u^\pm$ with $\norm{u^n-u^\pm}_p\le(1-\theta^p)^{-1/p}\,\theta^n\norm{u^0-Tu^0}_p$, where $\theta=\dfrac{2K\norm{u^0-Tu^0}_p+\kappa}{2K\norm{u^0-Tu^0}_p+1}$, and with asymptotic rate at most $\kappa$. If all $g_j$ are odd, then $u^-=-u^+$.
\end{mainthm}

\begin{proof}
By~\eqref{eq:Econd}, $R^p<\tfrac12(2a)^p$, so $R<2^{1-1/p}a<a$. For $u\in C_\pm$ we have $\abs{u_1\mp a}\le\norm{u\mp ae_1}_p\le R<a$, so $\sgn(u_1)=\pm1$ and $Tu=\pm ae_1+AG(u)$. By Lemma~\ref{lem:matrix}, $\norm{Tu\mp ae_1}_p\le\alpha\Gamma=R$, so $T(C_\pm)\subseteq C_\pm$. For $u,v$ in the same piece, $\norm{Tu-Tv}_p=\norm{A(G(u)-G(v))}_p\le\alpha L\norm{u-v}_p$. Now apply Theorem~\ref{thm:D} with $c_\pm=\pm ae_1$, $s=2a$ and $\lambda=\alpha L$; the conditions~\eqref{D1}--\eqref{D3} are exactly~\eqref{eq:Econd}. The estimate is Theorem~\ref{thm:A}(b) with $r=p$. If the $g_j$ are odd, then $T(-u)=-Tu$ and $-C_+=C_-$, so $-u^+$ is a fixed point in $C_-$, hence equal to $u^-$.
\end{proof}

\begin{example}\label{ex:lp}
Let $p=\tfrac12$ (so $K=2$), $a=8$, $\kappa=\tfrac1{32}$,
\[
a_{ij}=\tfrac1{16}\,4^{-(i-j)}\ (i\ge j),\qquad a_{ij}=0\ (i<j),\qquad g_j(t)=4^{-(j+1)}\sin t .
\]
Then $A$ is a discrete Volterra (lower-triangular) operator. We have $\sum_{i\ge j}\abs{a_{ij}}^{1/2}=\tfrac14\sum_{k\ge0}2^{-k}=\tfrac12$, so $\alpha=\tfrac14$. Also $L=\sup_j4^{-(j+1)}=\tfrac1{16}$ and $\Gamma=\left(\sum_j2^{-(j+1)}\right)^2=\tfrac14$, so $R=\tfrac1{16}$. The conditions~\eqref{eq:Econd} read
\[
\tfrac1{64}\,(1+16\cdot\tfrac1{16})=\tfrac1{32}\le\kappa,\qquad2\cdot\tfrac14=\tfrac12\le(1-2^{-1/2})\cdot4=4-2\sqrt2\approx1.17,\qquad16\ge8,
\]
and all hold. Since $a_{ij}\gamma_j=4^{-(i+3)}$, the system becomes
\[
u_i=8\,\sgn(u_1)\,\delta_{i1}+4^{-(i+3)}\sum_{j=1}^{i}\sin u_j,\qquad i\ge1 .
\]
By Theorem~\ref{thm:E} it has exactly two solutions $u^\pm$ in $C=\bar B(8e_1,\tfrac1{16})\cup\bar B(-8e_1,\tfrac1{16})$, and $u^-=-u^+$. Because the system is triangular, $u^+$ can be computed coordinate by coordinate:
\[
u^+_1\approx8.003862,\quad u^+_2\approx9.666\times10^{-4},\quad u^+_3\approx2.417\times10^{-4},\quad u^+_4\approx6.043\times10^{-5},\ \dots
\]
with $\norm{u^+-8e_1}_{1/2}\approx0.01546$ and $\norm{u^+-u^-}_{1/2}=2\norm{u^+}_{1/2}\approx16.72$. Starting from $u^0=8e_1$ we have $\norm{u^0-Tu^0}_{1/2}\approx0.01546$, $\theta\approx0.0877$ and $(1-\theta^{1/2})^{-2}\approx2.018$. Table~\ref{tab:lp} compares the a priori bound of Theorem~\ref{thm:A}(b) with the actual errors.

\begin{table}[ht]
\centering
\begin{tabular}{@{}ccc@{}}
\toprule
$n$ & $\norm{u^n-u^+}_{1/2}$ & a priori bound $(1-\theta^{1/2})^{-2}\theta^n\norm{u^0-Tu^0}_{1/2}$\\
\midrule
0 & $1.55\times10^{-2}$ & $3.12\times10^{-2}$\\
1 & $8.55\times10^{-6}$ & $2.73\times10^{-3}$\\
2 & $8.50\times10^{-9}$ & $2.40\times10^{-4}$\\
3 & $5.92\times10^{-12}$ & $2.10\times10^{-5}$\\
4 & $8.98\times10^{-15}$ & $1.84\times10^{-6}$\\
\bottomrule
\end{tabular}
\caption{Picard iteration for Example~\ref{ex:lp} (first 60 coordinates; the neglected tail is below $10^{-17}$).}
\label{tab:lp}
\end{table}

The actual convergence is faster than the bound because, within each well, $T$ has Lipschitz constant $\alpha L=1/64$. The bound of Theorem~\ref{thm:A} uses only~\eqref{eq:Rk}, and it holds uniformly over all starting points in $C$.
\end{example}

\begin{remark}[$N$-stable generalization]\label{rem:Nstable}
Theorem~\ref{thm:E} and Example~\ref{ex:lp} extend to an $N$-stable system by replacing the two-valued forcing $a\,\sgn(u_1)\delta_{i1}$ by a Voronoi-cell forcing on $N$ centers. Fix $N\ge 2$, scalars $\beta_1<\dots<\beta_N\in\R$ with $\gamma:=\min_{i\ne j}\abs{\beta_i-\beta_j}>0$, and $\sigma:\R\to\{1,\dots,N\}$ the piecewise-constant map $\sigma(t):=\arg\min_k\abs{t-\beta_k}$. The system $u_i=\beta_{\sigma(u_1)}\,\delta_{i1}+\sum_j a_{ij}\,g_j(u_j)$ then has exactly $N$ solutions in $C=\bigcup_i\bar B(\beta_i e_1,R)\subseteq\ell^p$, pairwise separated by $\gamma/2$, under the conditions $\alpha L(1+4^{1/p}R)\le\kappa$, $2R^p\le(1-2^{-p})\gamma^p$, $\gamma\ge 4(1+4^{1/p}R)$. The argument is Theorem~\ref{thm:E}'s proof, with $(C_+,C_-)$ replaced by $(C_i)_{i=1}^N$ and the sign condition $\sgn(u_1)=\pm 1$ for $u\in C_\pm$ replaced by $\sigma(u_1)=i$ for $u\in C_i$ (which follows from the bound $\abs{u_1-\beta_i}\le\norm{u-\beta_i e_1}_p$ (since $\abs{v_1}\le\norm v_p$ for every $v\in\ell^p$), together with $\norm{u-\beta_ie_1}_p\le R$ and $(2R)^p\le(1-2^{-p})\gamma^p\le 2^{-p}\gamma^p$ for $p\in(0,1]$, hence $R\le\gamma/4<\gamma/2$); then Corollary~\ref{cor:Npiece} gives the conclusion.
\end{remark}

\subsection{Volterra-type equations in \texorpdfstring{$H^p$, $0<p<1$}{H\textasciicircum p, 0<p<1}}\label{sec:Hp}

The second application takes place in the Hardy spaces $H^p$, $0<p<1$, with the Volterra-type operators $T_g$ as linear part. We perturb the linear equation $f=\pm a+T_gf$ by a nonlinear term that depends on the value of $f$ at a point, and show that the resulting equation again has exactly two solutions, which in a model case can be written down explicitly.

For $0<p<\infty$, $H^p$ is the space of analytic $f$ on the unit disc $\D$ with
\[
\norm f_p^p=\sup_{0\le r<1}\frac1{2\pi}\int_0^{2\pi}\abs{f(re^{it})}^p\,dt<\infty;
\]
see~\cite{Duren}. For $0<p<1$, $H^p$ is a $p$-Banach space, since $\abs{\alpha+\beta}^p\le\abs\alpha^p+\abs\beta^p$. For $g$ analytic on $\D$, put
\[
T_gf(z)=\int_0^zf(\zeta)\,g'(\zeta)\,d\zeta .
\]
This operator was introduced by Pommerenke~\cite{Pom}. Recall that $\BMOA$ is the space of analytic functions $g$ on $\D$ of bounded mean oscillation, that is, $g\in H^2$ and
\[
\norm g_{*}:=\sup_{a\in\D}\norm{g\circ\phi_a-g(a)}_{H^2}<\infty,\qquad \phi_a(z)=\frac{a-z}{1-\bar az};
\]
equivalently, the boundary function of $g$ belongs to $\mathrm{BMO}(\mathbb T)$, the John--Nirenberg space of bounded mean oscillation on the unit circle, or $\abs{g'(z)}^2(1-\abs z^2)\,dA(z)$ is a Carleson measure on $\D$ (that is, $\mu(S(I)):=\int_{S(I)}\abs{g'(z)}^2(1-\abs z^2)\,dA(z)\le C\,\abs I$ for every arc $I\subseteq\partial\D$, where $S(I)$ is the Carleson box over $I$). One has $H^\infty\subsetneq\BMOA$; for instance $\log\frac1{1-z}\in\BMOA\setminus H^\infty$. See~\cite{Garnett,Girela} for details. By Aleman--Siskakis~\cite{AS} ($1\le p<\infty$) and Aleman--Cima~\cite{AC} ($0<p<\infty$), $T_g$ is bounded on $H^p$ if and only if $g\in\BMOA$. For $g(z)=z$ it is the integration operator $J$, and for $g(z)=\log\frac1{1-z}$ it is closely related to the Ces\`aro operator. Note that $T_gf(0)=0$ and $T_{\varepsilon g}=\varepsilon T_g$.

\begin{lemma}\label{lem:pointeval}
For $0<p<\infty$, $f\in H^p$ and $z\in\D$, $\abs{f(z)}^p\le\norm f_p^p/(1-\abs z^2)$. In particular $\abs{f(0)}\le\norm f_p$.
\end{lemma}

\begin{proof}
Since $f$ is holomorphic, $\abs f^p$ is subharmonic on $\D$ (for a holomorphic $f$, $\log\abs f$ is subharmonic, hence so is $\abs f^p=e^{p\log\abs f}$), and subharmonicity gives $\abs{f(0)}^p\le M_p^p(r,f):=\frac1{2\pi}\int_0^{2\pi}\abs{f(re^{it})}^p\,dt\le\norm f_p^p$. For general $z$, let $\phi_z(w)=(z-w)/(1-\bar zw)$ and $F=(f\circ\phi_z)\,(\phi_z')^{1/p}$, using an analytic branch of $(\phi_z')^{1/p}$, which exists because $\phi_z'$ does not vanish on $\D$. Then $F\in H^p$ and $\norm F_p=\norm f_p$: on the circle, $\abs{F}^p=\abs{f\circ\phi_z}^p\abs{\phi_z'}$, and the substitution $\xi=\phi_z(e^{i\theta})$ has $\abs{d\xi}=\abs{\phi_z'(e^{i\theta})}\,d\theta$, so the boundary integrals of $\abs F^p$ and $\abs f^p$ agree (see~\cite[Ch.~2]{Duren}). Finally $\phi_z(0)=z$ and $\abs{\phi_z'(0)}=1-\abs z^2$, so $\abs{F(0)}=\abs{f(z)}(1-\abs z^2)^{1/p}\le\norm F_p=\norm f_p$.
\end{proof}

Let $g\in\BMOA$ with $\tau:=\norm{T_g}_{H^p\to H^p}<1$. Then $w:=\sum_{n\ge0}T_g^n1$ converges in $H^p$, because $\sum_n\norm{T_g^n1}_p^p\le\sum_n\tau^{np}<\infty$. It satisfies $w=1+T_gw$ and $w(0)=1$, so $\norm w_p\ge1$ by Lemma~\ref{lem:pointeval}. Let $z_0\in\D$, $k\in H^p$, and let $\psi:\C\to\C$ satisfy $\abs{\psi}\le\beta$ and $\abs{\psi(\xi)-\psi(\zeta)}\le\ell\abs{\xi-\zeta}$. For $a>0$ consider
\begin{equation}\label{eq:Hp}
f(z)=a\,\sgn\left(\mathrm{Re}\,f(0)\right)+\int_0^zf(\zeta)\,g'(\zeta)\,d\zeta+\psi\left(f(z_0)\right)\,k(z),\qquad z\in\D,
\end{equation}
that is, $f=Tf$ with $Tf:=\sigma(f)\,a+T_gf+\psi(f(z_0))\,k$ and $\sigma(f)=\sgn\mathrm{Re}\,f(0)$. The nonlinearity enters through a point evaluation, which is a continuous functional on $H^p$ by Lemma~\ref{lem:pointeval}.

\begin{mainthm}\label{thm:F}
Let $0<p<1$, $\kappa\in[0,1)$, $a,R>0$, and let $g,\tau,w,z_0,k,\psi,\beta,\ell$ be as above. Put
\[
\lambda:=\left(\tau^p+\frac{\ell^p\norm k_p^p}{1-\abs{z_0}^2}\right)^{1/p},
\]
and assume
\begin{align}
&\tau^pR^p+\beta^p\norm k_p^p\le R^p, \tag{F1}\label{F1}\\
&\lambda\,(1+4^{1/p}R)\le\kappa, \tag{F2}\label{F2}\\
&2R^p\le(1-2^{-p})(2a)^p, \tag{F3}\label{F3}\\
&2a\ge4(1+4^{1/p}R). \tag{F4}\label{F4}
\end{align}
Let $C_\pm=\left\lbrace f\in H^p:\norm{f\mp aw}_p\le R\right\rbrace $ and $C=C_+\cup C_-$. Then equation~\eqref{eq:Hp} has exactly two solutions $f_\pm$ in $C$, with $f_\pm\in C_\pm$ and $\norm{f_+-f_-}_p\ge a\norm w_p\ge a$. Picard iterates converge from every point of $C_\pm$ to $f_\pm$, with the estimates of Theorem~\ref{thm:A} for $b=2^{1/p-1}$ and $C_\theta=(1-\theta^p)^{-1/p}$. If $\psi$ is odd, then $f_-=-f_+$.
\end{mainthm}

\begin{proof}
As in the proof of Theorem~\ref{thm:E}, \eqref{F3} gives $R<a$. Let $f\in C_\pm$ and write $f=\pm aw+x$ with $\norm x_p\le R$. Then $f(0)=\pm a+x(0)$ and $\abs{x(0)}\le\norm x_p<a$, so $\sigma(f)=\pm1$. Using $w=1+T_gw$,
\[
Tf\mp aw=\pm a(1+T_gw-w)+T_gx+\psi(f(z_0))k=T_gx+\psi(f(z_0))\,k .
\]
By~\eqref{eq:pineq} and~\eqref{F1}, $\norm{Tf\mp aw}_p^p\le\tau^pR^p+\beta^p\norm k_p^p\le R^p$, so $T(C_\pm)\subseteq C_\pm$. For $f,h$ in the same piece,
\[
\norm{Tf-Th}_p^p\le\tau^p\norm{f-h}_p^p+\ell^p\abs{f(z_0)-h(z_0)}^p\norm k_p^p\le\lambda^p\norm{f-h}_p^p
\]
by Lemma~\ref{lem:pointeval}. Here $s=\norm{aw-(-aw)}_p=2a\norm w_p\ge2a$, and the conditions~\eqref{D2}--\eqref{D3} are increasing in $s$. So~\eqref{F3}--\eqref{F4} imply~\eqref{D2}--\eqref{D3}, while~\eqref{F2} is~\eqref{D1}. Theorem~\ref{thm:D} now gives all claims. For odd $\psi$, $T(-f)=-Tf$ and $-C_+=C_-$.
\end{proof}

For the integration operator $J=T_z$, $Jf(z)=\int_0^zf(\zeta)\,d\zeta$, the norm on $H^p$ can be bounded explicitly by an elementary argument.

\begin{lemma}\label{lem:J}
For every $0<p<\infty$,
\[
\norm J_{H^p\to H^p}\le\left(2(1+\sqrt2)^2\right)^{1/p}=(6+4\sqrt2)^{1/p}.
\]
\end{lemma}

\begin{proof}
Let $f\in H^p$ and let $f^+(e^{it})=\sup_{0\le\rho<1}\abs{f(\rho e^{it})}$ be its radial maximal function. For $0\le r<1$,
\[
\abs{Jf(re^{it})}=\left|\int_0^r f(\rho e^{it})\,e^{it}\,d\rho\right|\le r\,f^+(e^{it})\le f^+(e^{it}),
\]
so $\norm{Jf}_p\le\norm{f^+}_{L^p(\mathbb T)}$, where $\mathbb T$ carries the normalized measure $dt/2\pi$. By the Riesz factorization theorem~\cite[Ch.~2]{Duren}, every $f\in H^p$ ($0<p<\infty$) factors as $f=Bh$, where $B$ is a \emph{Blaschke product} (an analytic function on $\D$ of the form $B(z)=z^{m}\prod_n\tfrac{\bar z_n}{\abs{z_n}}\tfrac{z_n-z}{1-\bar z_nz}$ with $\abs B\le 1$ on $\D$ and $\abs B=1$ a.e.\ on $\partial\D$) and $h\in H^p$ is zero-free, with $\norm h_p=\norm f_p$. Then $G:=h^{p/2}$ is a well-defined element of $H^2$ (using a holomorphic branch of $h^{p/2}$, available because $h$ is zero-free and $\D$ is simply connected) with $\norm G_2^2=\norm h_p^p=\norm f_p^p$, and $f^+\le h^+=(G^+)^{2/p}$ because $\abs B\le1$, so that $\norm{f^+}_p^p\le\norm{G^+}_2^2$. Let $G^*\in L^2(\mathbb T)$ be the boundary function of $G$. Since $G$ is the Poisson integral of $G^*$ and the Poisson kernel is even, positive, decreasing in $\abs t$ on $[0,\pi]$ and of mass one, $G^+\le M_{\mathbb T}G^*$, where $M_{\mathbb T}G^*(x):=\sup_{0<h\le\pi}\tfrac1{2h}\int_{x-h}^{x+h}\abs{G^*(e^{it})}\,dt$ is the centered Hardy--Littlewood maximal operator on $\mathbb T$~\cite[Ch.~1]{Duren}; also $\norm{G^*}_{L^2}=\norm G_2$. Identify $\mathbb T$ with $[-\pi,\pi)$ and extend $G^*$ periodically. For $x\in[-\pi,\pi)$, every arc used in $M_{\mathbb T}G^*(x)$ is an interval $[x-h,x+h]\subseteq[-2\pi,2\pi]$, and the averages over it are the same on $\mathbb T$ and on $\R$, so $M_{\mathbb T}G^*(x)\le M_{\R}(G^*\mathbf 1_{[-2\pi,2\pi]})(x)$, where $M_\R$ is the uncentered maximal operator on $\R$. By Grafakos--Montgomery-Smith~\cite{GMS}, $\norm{M_\R}_{L^2(\R)\to L^2(\R)}=1+\sqrt2$. Hence
\[
\int_{-\pi}^{\pi}(M_{\mathbb T}G^*)^2\,dt\le(1+\sqrt2)^2\int_{-2\pi}^{2\pi}\abs{G^*}^2\,dt=2(1+\sqrt2)^2\int_{-\pi}^{\pi}\abs{G^*}^2\,dt,
\]
that is, $\norm{G^+}_2^2\le2(1+\sqrt2)^2\norm G_2^2$. Combining, $\norm{Jf}_p^p\le2(1+\sqrt2)^2\norm f_p^p$.
\end{proof}

\begin{remark}\label{rem:bmoa-bv}
The constant in Lemma~\ref{lem:J} is far from sharp; we only need it to be explicit, because the hypotheses of Theorem~\ref{thm:F} contain $\tau=\norm{T_g}_{H^p\to H^p}$ and must be checkable in applications. By~\cite{AC}, $T_g$ is bounded on $H^p$ for every $g\in\BMOA$, but the available constants are not explicit in $p$. The same argument as in Lemma~\ref{lem:J}, using $\abs{T_gf(re^{it})}\le f^+(e^{it})\int_0^1\abs{g'(\rho e^{it})}\,d\rho$, gives
\[
\norm{T_g}_{H^p\to H^p}\le(6+4\sqrt2)^{1/p}\,V(g),\qquad V(g):=\sup_{t}\int_0^1\abs{g'(\rho e^{it})}\,d\rho,
\]
which supplies an explicit upper bound on $\tau$ whenever $g$ has bounded radial variation. The class $\{g:V(g)<\infty\}$ is a strict subclass of $\BMOA$ (for instance $\log\tfrac1{1-z}\in\BMOA$ has $V(g)=\infty$) that includes the model symbol $g(z)=\varepsilon z$ of Example~\ref{ex:Hp} below. In this sense Theorem~\ref{thm:F} is quantitative for symbols of bounded radial variation and qualitative for general $g\in\BMOA$: the existence, multiplicity and separation statements hold in the full $\BMOA$ setting whenever the hypotheses~\eqref{F1}--\eqref{F4} can be verified, but checking them requires an explicit upper bound on $\tau$, which is supplied here only in the bounded radial variation case.
\end{remark}

\begin{remark}
Since $T_{\varepsilon g}=\varepsilon T_g$, for any $g\in\BMOA$, $z_0\in\D$ and $k\in H^p$, conditions \eqref{F1}--\eqref{F4} hold for all sufficiently small $\abs\varepsilon$ and $\beta,\ell$ once $R$ and $a$ are fixed accordingly. In general the only non-explicit quantity is $\norm{T_g}_{H^p\to H^p}$, which by~\cite{AC} is comparable to the $\BMOA$ seminorm of $g$ with constants depending on $p$; for $g(z)=\varepsilon z$ it is bounded explicitly by Lemma~\ref{lem:J}.
\end{remark}

\begin{example}\label{ex:Hp}
Let $p=\tfrac12$, $g(z)=\varepsilon z$ (so $T_g=\varepsilon J$), $z_0=\tfrac12$, $k(z)=z$ (so $\norm k_{1/2}=1$), $\psi(\zeta)=\eta\sin(\mathrm{Re}\,\zeta)$ with $\eta=\tfrac1{100}$ (so $\beta=\ell=\eta$ and $\psi$ is odd), $a=8$, $R=\tfrac1{16}$, $\kappa=\tfrac14$. Take $\varepsilon=2\cdot10^{-4}$. By Lemma~\ref{lem:J}, $\norm J_{H^{1/2}\to H^{1/2}}\le(6+4\sqrt2)^2=68+48\sqrt2<136$, so $\tau\le\varepsilon\cdot136<\tfrac1{25}$ and $\tau^{1/2}\le\tfrac15$. Then:
\begin{itemize}
\item \eqref{F1}: $(\tau R)^{1/2}+\eta^{1/2}\le\tfrac1{20}+\tfrac1{10}=\tfrac3{20}\le\tfrac14=R^{1/2}$;
\item \eqref{F2}: $\lambda^{1/2}=\tau^{1/2}+\tfrac43\eta^{1/2}\le\tfrac15+\tfrac2{15}=\tfrac13$, so $\lambda\le\tfrac19$ and $\lambda(1+16R)\le\tfrac29\le\tfrac14$;
\item \eqref{F3}--\eqref{F4}: as in Example~\ref{ex:lp}.
\end{itemize}
Here $w=e^{\varepsilon z}$. Equation~\eqref{eq:Hp} reads
\[
f(z)=8\,\sgn(\mathrm{Re}\,f(0))+\varepsilon\int_0^zf(\zeta)\,d\zeta+\tfrac1{100}\sin\left(\mathrm{Re}\,f(\tfrac12)\right)\,z .
\]
By Theorem~\ref{thm:F} it has exactly two solutions in $C=\left\lbrace \norm{f-8e^{\varepsilon z}}_{1/2}\le\tfrac1{16}\right\rbrace \cup\left\lbrace \norm{f+8e^{\varepsilon z}}_{1/2}\le\tfrac1{16}\right\rbrace $, and $f_-=-f_+$. They can be found explicitly. Put $c:=\eta\sin(\mathrm{Re}\,f_+(\tfrac12))\in\R$. Differentiating gives $f_+'=\varepsilon f_++c$ with $f_+(0)=8$, hence
\[
f_+(z)=\left(8+\frac c\varepsilon\right)e^{\varepsilon z}-\frac c\varepsilon,\qquad c=\eta\sin\left(\left(8+\frac c\varepsilon\right)e^{\varepsilon/2}-\frac c\varepsilon\right).
\]
Theorem~\ref{thm:F} guarantees that exactly one root $c$ of this scalar equation gives $f_+\in C_+$. Numerically, for $\varepsilon=2\cdot10^{-4}$ the iteration $c\mapsto\eta\sin\left((8+c/\varepsilon)e^{\varepsilon/2}-c/\varepsilon\right)$ started at $c=0$ converges to $c\approx0.0098851$, and $\mathrm{Re}\,f_+(\tfrac12)\approx8.00574$. As $\varepsilon\to0$, $f_+(z)\to8+c_0z$ with $c_0=\eta\sin(8+c_0/2)\approx0.0098863$. As a check, $f_+-8e^{\varepsilon z}=c\,(e^{\varepsilon z}-1)/\varepsilon=c\sum_{n\ge1}\varepsilon^{n-1}z^n/n!$, and since $\norm{z^n}_{1/2}=1$, the $\tfrac12$-subadditivity of $\norm\cdot_{1/2}^{1/2}$ gives $\norm{f_+-8e^{\varepsilon z}}_{1/2}^{1/2}\le\abs c^{1/2}\sum_{n\ge1}\left(\varepsilon^{n-1}/n!\right)^{1/2}\le1.011\,\abs c^{1/2}$. Hence $\norm{f_+-8e^{\varepsilon z}}_{1/2}\le1.03\,\abs c<0.011<R$, so indeed $f_+\in C_+$.
\end{example}

\section{Concluding remarks}\label{sec:open}

Theorems~\ref{thm:KAC}--\ref{thm:old23} establish existence of fixed points for rational-type contractions. The present paper completes them as follows. A map satisfying~\eqref{eq:Rk} is automatically continuous, its Picard iterates converge at an explicit rate, its fixed points are uniformly separated and each attracts a ball of fixed radius, and its basins of attraction are open and closed. Consequently the fixed point is unique on connected spaces, and several fixed points can occur only on disconnected spaces, where they do occur by the gluing principle and the two applications. The quasi-Banach setting differs from the metric one (Example~\ref{ex:incomparable}).

The natural next step is the set-valued case. For $F$ mapping $U$ to its nonempty compact subsets one considers
\[
H(Fu,Fv)\le\frac{D(u,Fv)+D(v,Fu)+\kappa}{b\,[D(u,Fu)+D(v,Fv)]+1}\,d(u,v),
\]
where $H$ is the Pompeiu--Hausdorff $b$-metric and $D(x,A)=\min_{a\in A}d(x,a)$. The structure result of Theorem~\ref{thm:B}(f) does not carry over to fixed points $u\in Fu$: on $U=[0,1]$ the constant map $Fu=[0,1]$ satisfies this inequality, since $H(Fu,Fv)=0$, and every point of the connected space $U$ is a fixed point. The appropriate objects are the endpoints of $F$, that is, the points with $Fu=\{u\}$. Convergence rates, stability and the structure of the set of endpoints for such maps will be studied in a forthcoming paper.

\subsection*{Competing interests} The authors declare that they have no competing interests.


\begin{thebibliography}{99}

\bibitem{AC}
A.~Aleman and J.~A.~Cima, An integral operator on $H^p$ and Hardy's inequality, \emph{J. Anal. Math.} \textbf{85} (2001), 157--176.

\bibitem{AS}
A.~Aleman and A.~G.~Siskakis, An integral operator on $H^p$, \emph{Complex Variables Theory Appl.} \textbf{28} (1995), 149--158.

\bibitem{AGAY}
A.~H.~Ansari Komachali, Y.~U.~Gaba, M.~Aphane and I.~Yildirim, Refinement on fixed point results in metric and $b$-metric spaces, \emph{Filomat} \textbf{37}:22 (2023), 7581--7588. \url{https://doi.org/10.2298/FIL2322581A}

\bibitem{Aoki}
T.~Aoki, Locally bounded linear topological spaces, \emph{Proc. Imp. Acad. Tokyo} \textbf{18} (1942), 588--594.

\bibitem{Bakhtin}
I.~A.~Bakhtin, The contraction mapping principle in quasimetric spaces (Russian), \emph{Func. An., Gos. Ped. Inst. Unianowsk} \textbf{30} (1989), 26--37.

\bibitem{BW}
D.~W.~Boyd and J.~S.~W.~Wong, On nonlinear contractions, \emph{Proc. Amer. Math. Soc.} \textbf{20} (1969), 458--464.

\bibitem{Czerwik93}
S.~Czerwik, Contraction mappings in $b$-metric spaces, \emph{Acta Math. Inform. Univ. Ostraviensis} \textbf{1} (1993), 5--11.

\bibitem{Czerwik98}
S.~Czerwik, Nonlinear set-valued contraction mappings in $b$-metric spaces, \emph{Atti Sem. Mat. Fis. Univ. Modena} \textbf{46} (1998), 263--276.

\bibitem{Day}
M.~M.~Day, The spaces $L^p$ with $0<p<1$, \emph{Bull. Amer. Math. Soc.} \textbf{46} (1940), 816--823.

\bibitem{DV}
M.~Demma and P.~Vetro, Picard sequence and fixed point results on $b$-metric spaces, \emph{J. Funct. Spaces} \textbf{2015} (2015), Art. ID 189861.

\bibitem{Duren}
P.~L.~Duren, \emph{Theory of $H^p$ Spaces}, Academic Press, New York, 1970.

\bibitem{Garnett}
J.~B.~Garnett, \emph{Bounded Analytic Functions}, revised 1st ed., Graduate Texts in Mathematics 236, Springer, New York, 2007.

\bibitem{Girela}
D.~Girela, Analytic functions of bounded mean oscillation, in: \emph{Complex Function Spaces (Mekrij\"arvi, 1999)}, Univ. Joensuu Dept. Math. Rep. Ser. 4, Univ. Joensuu, Joensuu, 2001, pp.~61--170.

\bibitem{GMS}
L.~Grafakos and S.~Montgomery-Smith, Best constants for uncentred maximal functions, \emph{Bull. London Math. Soc.} \textbf{29} (1997), 60--64.

\bibitem{HDR}
H.~Huang, G.~Deng and S.~Radenovi\'c, Fixed point theorems in $b$-metric spaces with applications to differential equations, \emph{J. Fixed Point Theory Appl.} \textbf{20} (2018), Paper No.~52.

\bibitem{JKR}
M.~Jovanovi\'c, Z.~Kadelburg and S.~Radenovi\'c, Common fixed point results in metric-type spaces, \emph{Fixed Point Theory Appl.} \textbf{2010} (2010), Art. ID 978121.

\bibitem{KPR}
N.~J.~Kalton, N.~T.~Peck and J.~W.~Roberts, \emph{An $F$-Space Sampler}, London Math. Soc. Lecture Note Ser. 89, Cambridge Univ. Press, 1984.

\bibitem{KAC}
F.~Khojasteh, M.~Abbas and S.~Costache, Two new types of fixed point theorems in complete metric spaces, \emph{Abstr. Appl. Anal.} \textbf{2014} (2014), Art. ID 325840.

\bibitem{KShR}
F.~Khojasteh, S.~Shukla and S.~Radenovi\'c, A new approach to the study of fixed point theory for simulation functions, \emph{Filomat} \textbf{29}:6 (2015), 1189--1194.

\bibitem{KS}
W.~Kirk and N.~Shahzad, \emph{Fixed Point Theory in Distance Spaces}, Springer, Cham, 2014.

\bibitem{MK}
A.~Meir and E.~Keeler, A theorem on contraction mappings, \emph{J. Math. Anal. Appl.} \textbf{28} (1969), 326--329.

\bibitem{MM}
R.~Miculescu and A.~Mihail, New fixed point theorems for set-valued contractions in $b$-metric spaces, \emph{J. Fixed Point Theory Appl.} \textbf{19} (2017), 2153--2163.

\bibitem{PS}
M.~Paluszy\'nski and K.~Stempak, On quasi-metric and metric spaces, \emph{Proc. Amer. Math. Soc.} \textbf{137} (2009), 4307--4312.

\bibitem{Pom}
Ch.~Pommerenke, Schlichte Funktionen und analytische Funktionen von beschr\"ankter mittlerer Oszillation, \emph{Comment. Math. Helv.} \textbf{52} (1977), 591--602.

\bibitem{Rolewicz}
S.~Rolewicz, On a certain class of linear metric spaces, \emph{Bull. Acad. Polon. Sci. Cl. III} \textbf{5} (1957), 471--473.

\bibitem{Rus}
I.~A.~Rus, Remarks on Ulam stability of the operatorial equations, \emph{Fixed Point Theory} \textbf{10} (2009), 305--320.

\bibitem{Suzuki}
T.~Suzuki, Basic inequality on a $b$-metric space and its applications, \emph{J. Inequal. Appl.} \textbf{2017} (2017), Paper No.~256.

\bibitem{Wardowski}
D.~Wardowski, Fixed points of a new type of contractive mappings in complete metric spaces, \emph{Fixed Point Theory Appl.} \textbf{2012} (2012), Paper No.~94.

\end{thebibliography}
\end{document}